\documentclass[12pt]{article}
\usepackage{amsmath,amssymb}

\newtheorem{theorem}{Theorem}

\newtheorem{lemma}[theorem]{Lemma}

\newtheorem{proposition}[theorem]{Proposition}

\newtheorem{remark}[theorem]{Remark}

\newdimen\dummy
\dummy=\oddsidemargin
\let\text=\mbox

\input{tcilatex}

\begin{document}

\title{Adaptive Nonparametric Regression on Spin Fiber Bundles}
\author{Claudio Durastanti \\
%EndAName
\textit{Dipartimento di Matematica, Universit\`{a} di Pavia} \and Daryl
Geller \\
%EndAName
\textit{Department of Mathematics, Stony Brook University} \and Domenico
Marinucci \\
%EndAName
\textit{Dipartimento di Matematica, Universit\`{a} di Roma Tor Vergata} }
\maketitle

\begin{abstract}
The construction of adaptive nonparametric procedures by means of wavelet
thresholding techniques is now a classical topic in modern mathematical
statistics. In this paper, we extend this framework to the analysis of
nonparametric regression on sections of spin fiber bundles defined on the
sphere. This can be viewed as a regression problem where the function to be
estimated takes as its values algebraic curves (for instance, ellipses)
rather than scalars, as usual. The problem is motivated by many important
astrophysical applications, concerning for instance the analysis of the weak
gravitational lensing effect, i.e. the distortion effect of gravity on the
images of distant galaxies. We propose a thresholding procedure based upon
the (mixed) spin needlets construction recently advocated by Geller and
Marinucci (2008,2010) and Geller et al. (2008,2009), and we investigate
their rates of convergence and their adaptive properties over spin Besov
balls.
\end{abstract}

\medskip


\medskip


\section{Introduction}

\label{sec:introduction}

Over the last two decades, wavelet techniques have become a well-established
tool for the analysis of statistical nonparametric problems, especially in
the framework of minimax estimation. The seminal contribution in this area
was provided by Donoho et al. in \cite{donoho1}, where it was proved that
nonlinear wavelet estimators based on thresholding techniques achieve nearly
optimal minimax rates (up to logarithmic terms) for a wide class of
nonparametric estimation of unknown density and regression functions. The
theory has been enormously developed ever since - we refer to \cite{WASA}
for a textbook reference.\newline
The bulk of this literature has focussed on estimation in standard Euclidean
frameworks, such as $\mathbb{R}$ or $\mathbb{R}^{n}$. More recently,
applications from various scientific fields have drawn a lot of attention on
more general settings, such as spherical data or more general manifolds (see %
\cite{randomfields}). This environment has recently experienced a remarkable
amount of activity, both from the purely mathematical point of view and in
terms of applications to empirical data. \newline
In particular, a highly successful construction of a second-generation
wavelet system on the sphere (the so-called needlets) has been introduced by %
\cite{npw1}, \cite{npw2}; this approach has been extended to more general
manifolds and unbounded support in the harmonic domain by \cite{gm1}, \cite%
{gm2}, \cite{gm3}. The investigation of the stochastic properties of
needlets when implemented on spherical random fields is due to \cite{bkmpAoS}%
, \cite{bkmpBer}, \cite{ejslan}, \cite{spalan} and \cite{mayeli}, where
applications to several statistical procedures are also considered. These
procedures have been mainly motivated by issues arising in Cosmology and
Astrophysics, and indeed several applications to experimental data have
already been implemented: for instance, those from the satellite WMAP
mission from NASA, focussing on the so-called Cosmic Microwave Background
radiation, see \cite{pbm06}, \cite{mpbb08}, \cite{pietrobon1}, \cite{fay08}, %
\cite{pietrobon2}, \cite{rudjord1}, \cite{dela08}, \cite{rudjord2}, \cite%
{cabella10} and \cite{ghosh}. These applications, however, have not been
focussed on thresholding estimates and minimax results, but rather to random
fields issues, such as angular power spectrum estimation, higher-order
spectra, testing for Gaussianity and isotropy, and several others (see also %
\cite{m},\cite{cama}).\newline
More recently, a few papers have focussed on the use of needlets to develop
estimators within the thresholding paradigm, in the framework of directional
data. The pioneering contribution here is due to \cite{bkmpAoSb}, see also %
\cite{thres}, \cite{Kerkypicard}, \cite{knp}; applications to astrophysical
data is still under way. Earlier results on minimax estimators for spherical
data, outside the needlets approach, are due to Kim and coauthors (see \cite%
{kim}, \cite{kimkoo}, \cite{kookim}).\newline
Another important generalization of the needlet approach has been recently
advocated by \cite{gelmar}; applications to statistics can be found in \cite%
{glm}. This development is again motivated by Cosmology and Astrophysics. In
particular, we noted above as some extremely influential satellite missions
from NASA and ESA (WMAP and Planck, respectively) are currently collecting
data on the so-called Cosmic Microwave Background radiation, which can be
viewed as the realization of a scalar, isotropic, mean-square continuous
spherical random field (see for instance \cite{dode2004} for a review).
These same experiments are also collecting data on a much more elusive
cosmological feature, the so-called polarization of CMB. The latter can be
loosely described as observations on random ellipses living on the tangent
planes for each location on the celestial sphere. Mathematically, this can
be expressed by defining random sections of spin fiber bundles, a
generalization of the notion of scalar random fields (see \cite{gelmar}, %
\cite{glm}, \cite{gelmar2010} and Sections \ref{sec:spinfunctions}
below for much more details and discussion). Quite interestingly,
exactly the same mathematical framework describes the so-called weak
gravitational lensing induced on the observed shape of distant
galaxies by clusters of matter. This is again a major issue in the
analysis of astrophysical data (see for instance
\cite{bridle,great10} and the references therein): huge amount of
observational data are expected in the next decade, by means of
satellite missions in preparations such as Euclid.\newline
The applications of spin needlets to CMB polarization data is discussed in %
\cite{ghmkp}. The characterization of spin Besov spaces by means of needlets
decompositions is discussed by \cite{bkmp09} and \cite{gelmar2010}; the
latter reference also introduces an alternative construction for needlets on
spin fiber bundles (so-called mixed needlets), and provide its analytical
and statistical properties. \newline
Our purpose in this article is to exploit these results and classical
techniques to introduce and develop spin nonparametric regression, with a
view to applications to polarization and weak lensing data. In particular,
we investigate the properties of nonlinear hard thresholding estimates, and
we establish rates of convergence over a wide class of $L_{s}^{p}$ norms and
spin Besov spaces (see again \cite{bkmp09}, \cite{gelmar2010} and the
sections to follow for more detailed definitions). More precisely, we shall
assume to have observations on independent pairs of random variables,
respectively scalar and spin, $\left( X_{i},Y_{i;s}\right) $, $i=1,\ldots ,n$%
, $\left( X_{i}\right) \in \mathbb{S}^{2}$; we view $\left( X_{i}\right) $
as uniform random locations on the sphere, which correspond for instance to
the positions of observed galaxies. We shall then be concerned with the
regression model:%
\begin{equation}
Y_{i;s}=F_{s}\left( X_{i}\right) +\varepsilon _{i;s}\text{ ,}
\label{eqn:regression}
\end{equation}%
where $F_{s}\left( \cdot \right) $ is an unknown section of a spin fiber
bundle; for instance, for $s=2$ $F_{s}$ can be taken to represent the
geometric effect of the gravitational shear. We assume that this section
belongs to $L_{s}^{p}\left( \mathbb{S}^{2}\right) $, the space of the spin $%
s,$ $p$-integrable sections on the sphere. On the other hand, we assume the $%
\varepsilon _{i;s}$ are i.i.d. spin random variables, which can be viewed as
an observational error (to be interpreted, for instance, as the intrinsic
shape of the galaxy). We are then led to nonparametric estimation over an
unknown functional class, and we aim at procedures which are robust (i.e.
nearly optimal) for a wide class of $L_{s}^{p}$ norms, $1\leq p\leq \infty $%
. To address this issue, and given the properties of (mixed) spin needlets
established in \cite{gelmar}, \cite{gelmar2010}, we follow a classical
approach, as discussed for the classical case on $\mathbb{R}$ by Donoho et
al. (\cite{donoho1}), Hardle et al. (\cite{WASA}), (\cite{thres}), and many
other papers, see for instance (\cite{brown},\cite{cai}, \cite{tsyb}) for
some recent developments . In particular, as mentioned before we introduce
thresholding estimates and establish convergence rates for the resulting
nonlinear estimators. We stress that we consider at the same time estimators
based upon both spin constructions we have mentioned before, i.e. pure and
mixed spin needlets; the results with the two approaches are identical.
Sharp adaptation results for nonparametric regression on vector bundles were
recently established in an important paper by \cite{kim}. These authors
focus on the $p=2,$ and therefore exploit Fourier methods rather than
wavelets thresholding. For $s=1,$ our method can be viewed as a form of
adaptive regression for vector fields, and in this sense it relates also to
recent work on filament estimation by \cite{gpvw1}, \cite{gpvw2}. See also %
\cite{schwartzman} for some recent work on statistical analysis for
tensor-valued data.\newline
The plan of the paper is as follows. In Section 2 we review some background
material on spin fiber bundles, while in Section 3 we recall the
construction of spin and mixed spin needlets; for both sections we follow
closely earlier references, in particular \cite{gelmar} and \cite{gelmar2010}%
. In Section 4 we review some crucial material on spin Besov spaces, as
discussed earlier by \cite{bkmp09} and \cite{gelmar2010}. Section 5 and 6
include the most important contributions of this paper, namely the
presentation of the thresholding procedure and the investigation of its
asymptotic properties.

\section{Spin functions}

\label{sec:spinfunctions}

\subsection{Background and definitions}

The purpose of this Section is to review some background material on spin
fiber bundles; our presentation follows closely \cite{gelmar}, \cite{glm},%
\cite{gelmar2010}, to which we refer for more discussion and details.\newline
The concept of a spin function was introduced in the sixties by Newman and
Penrose in \cite{newman}, while working on gravitational radiation, see also %
\cite{goldb}, \cite{edth}. Writing in a physicists' jargon, they said that a
function $\eta $ has an integer-valued spin weight $s$ (or, briefly, that $%
\eta $ is a spin $s$ quantity) if, whenever a tangent vector at point $x\in
\mathbb{S}^{2}$ is rotated by an angle $\psi $ under a coordinate change, $%
\eta $ transforms as $\eta ^{\prime }=e^{is\psi }\eta $. This same idea is
formalized as follows by \cite{gelmar}). Let $U_{I}:=\mathbb{S}^{2}/\left\{
N,S\right\} $ be the chart that covers the sphere with the North and the
South poles subtracted: here we adopt the usual angular coordinates $\left(
\vartheta ,\varphi \right) $, $\vartheta \in \left( 0,\pi \right) $ and $%
\varphi \in \left[ -\pi ,\pi \right] $. Define the rotated charts $%
U_{R}=RU_{I}$, where $R\in SO\left( 3\right) $ (the special group of
rotations) and label the corresponding coordinates $\left( \vartheta
_{R},\varphi _{R}\right) $. For any $x\in \mathbb{S}^{2}$, we can fix a
''reference direction'' in the tangent plane at $x$ (labelled as usual $%
T_{x}\left( \mathbb{S}^{2}\right) $) by considering $\rho _{I}\left(
x\right) =\partial /\partial \varphi $, the unitary tangent vector in the
direction of the circle where $\vartheta $ is constant and $\varphi $ is
increasing. For every $x$ belonging to the intersection between the charts
corresponding to $U_{R}$ and $U_{I}$, we can uniquely measure the angle
associated to a change of coordinate by considering the angle between the
reference vector in the map $U_{I}$, and the reference vector in the rotated
chart, namely $\rho _{R}\left( x\right) =\partial /\partial \varphi _{R}$.
More generally, given $x\in \mathbb{S}^{2}$ and two charts $U_{R_{1}}$ and $%
U_{R_{2}}$ such that $x\in U_{R_{1}}\cap U_{R_{2}}$, the angle between $%
U_{R_{1}}$ and $U_{R_{2}}$, $\psi _{x,R_{1},R_{2}}$ is defined as the angle
between $\rho _{R_{1}}\left( x\right) $ and $\rho _{R_{2}}\left( x\right) $,
see \cite{gelmar}, \cite{glm} for a discussion on the orientation of this
angle. \newline
Fix now an open subset $G\subset \mathbb{S}^{2}$. The collection of
functions $\left\{ F_{R}\right\} _{R\in SO\left( 3\right) }$ is a spin $s$
function $F_{s}$ if and only if $\forall R_{1},R_{2}\in SO\left( 3\right) $
and all $x\in U_{R_{1}}\cap U_{R_{2}}\cap G$ we have:%
\begin{equation*}
F_{R_{2}}=e^{is\psi _{x,R_{1},R_{2}}}F_{R_{1}}
\end{equation*}%
We write $F_{s}\in C_{s}^{\infty }\left( G\right) $, if for every $R\in
SO\left( 3\right) $ the application $x\rightarrow F_{s}\left( x\right) $ is
smooth. Note that for $s=0$ we are back to the usual scalar functions.
\newline
From a differential geometry point of view, $C_{s}^{\infty }$ is the space
of sections over $G$ of the complex line bundle over the sphere $\mathbb{S}%
^{2}$ (see also \cite{leosa}, \cite{mal} for more discussion on this point
of view). The functional spaces $L_{s}^{p}\left( \mathbb{S}^{2}\right) $ are
then defined as%
\begin{equation*}
F_{s}\in L_{s}^{p}\left( \mathbb{S}^{2}\right) \Leftrightarrow \left\|
F_{s}\right\| _{L_{s}^{p}\left( \mathbb{S}^{2}\right) }=\left( \int_{\mathbb{%
S}^{2}}\left| F_{s}(x)\right| ^{p}dx\right) ^{1/p}<\infty
\end{equation*}%
Note that, while $F_{s}\left( x\right) $ is a section of the fiber bundle on
$\mathbb{S}^{2}$, $\left| F_{s}\left( x\right) \right| $ is a real valued
function on the sphere, because the modulus of $F_{s}$ does not depend on
the choice of the coordinate system: therefore the $L_{s}^{p}\left( \mathbb{S%
}^{2}\right) $ is well defined.

\subsection{Spin Spherical Harmonics}

\label{sec:spinharmonics}

We start by recalling the well-known expression for the spherical Laplacian $%
\Delta _{\mathbb{S}^{2}},$
\begin{equation*}
\Delta _{\mathbb{S}^{2}}:=\frac{1}{\sin ^{2}\vartheta }\frac{\partial }{%
\partial \varphi ^{2}}+\frac{1}{\sin \vartheta }\frac{\partial }{\partial
\vartheta }\left\{ \sin \vartheta \frac{\partial }{\partial \vartheta }%
\right\} \text{ .}
\end{equation*}%
A complete orthonormal set of eigenfunctions for the spherical Laplacian is
provided by the family of spherical harmonics $\left\{ Y_{lm}\right\} $, $%
l=0,1,2,...,m=-l,...,l:$%
\begin{equation*}
\Delta _{\mathbb{S}^{2}}Y_{lm}=-l(l+1)Y_{lm}\text{ , }\int_{\mathbb{S}%
^{2}}Y_{lm}(x)\overline{Y}_{lm}(x)dx=\delta _{l}^{l\prime }\delta
_{m}^{m\prime }\text{ ,}
\end{equation*}%
where $\delta _{a}^{b}$ denotes the Kronecker delta function. In the
spherical coordinates $(\vartheta ,\varphi )$
\begin{eqnarray*}
Y_{lm}(\vartheta ,\varphi ) &=&e^{im\varphi }\sqrt{\frac{2l+1}{4\pi }\frac{%
(l-m)!}{(l+m)!}}P_{lm}(\cos \vartheta )\text{ ,} \\
P_{lm}(x) &=&(1-x^{2})^{m/2}\frac{d}{dx^{m}}P_{l}(x)\text{ ,}
\end{eqnarray*}%
where $P_{l}(x)$ denotes the Legendre polynomials, see for instance (\cite%
{vmk}) for more analytic expressions and discussion. Denoting by $\left\{
\mathcal{H}_{l}\right\} $ the linear spaces spanned by the spherical
harmonics, the following decomposition holds (see for instance \cite%
{randomfields}):
\begin{equation*}
L^{2}\left( \mathbb{S}^{2}\right) =\bigoplus_{l\geq 0}{\mathcal{H}_{l}}\text{
,}
\end{equation*}%
that is, in the $L^{2}$ sense, for all $f\in L^{2}\left( \mathbb{S}%
^{2}\right) $%
\begin{equation*}
f\left( x\right) =\sum_{l,m}a_{lm}Y_{lm}\left( x\right) \text{ , }%
a_{lm}=\int_{\mathbb{S}^{2}}f(x)\overline{Y}_{lm}(x)dx\text{{\ .}}
\end{equation*}%
It is possible to introduce spin spherical harmonics as the eigenfunctions
of a second-order differential operator which generalizes the spherical
Laplacian (refer again to \cite{wiaux06},\cite{gelmar} for more details). To
this aim, consider the (\emph{spin raising} and \emph{spin lowering})
operators $\eth $ and $\overline{\eth }$, whose action on a spin function $%
F_{s}\left( \cdot \right) $ is provided by:
\begin{eqnarray*}
&&\eth F_{s}\left( \vartheta ,\varphi \right) =-\left( \sin \left( \theta
\right) \right) ^{s}\left[ \frac{\partial }{\partial \vartheta }+\frac{i}{%
\sin \left( \vartheta \right) }\frac{\partial }{\partial \varphi }\right]
\left( \sin \left( \theta \right) \right) ^{-s}F_{s}\left( \vartheta
,\varphi \right) \text{ ,} \\
&&\overline{\eth }F_{s}\left( \vartheta ,\varphi \right) =-\left( \sin
\left( \theta \right) \right) ^{-s}\left[ \frac{\partial }{\partial
\vartheta }-\frac{i}{\sin \left( \vartheta \right) }\frac{\partial }{%
\partial \varphi }\right] \left( \sin \left( \theta \right) \right)
^{s}F_{s}\left( \vartheta ,\varphi \right) \text{ .}
\end{eqnarray*}%
It should be noted that $\eth $ transforms spin $s$ functions into spin $s+1$
functions, $\eth C_{s}^{\infty }\rightarrow C_{s+1}^{\infty },$ while $%
\overline{\eth }$ transforms spin $s$ functions into spin $s-1$ functions, $%
\overline{\eth }C_{s}^{\infty }\rightarrow C_{s-1}^{\infty },$ which
justifies their names. The previous expressions should be written more
rigorously in terms of $\eth _{R},\overline{\eth }_{R},\vartheta
_{R},\varphi _{R},F_{s;R}$, because both the operators and the spin
functions depend on the choice of coordinates. More important, $\eth ,%
\overline{\eth }$ can be used to define a differential operator $\eth
\overline{\eth }$, which can be viewed as a generalization of the scalar
spherical Laplacian; indeed%
\begin{equation*}
-\overline{\eth }\eth Y_{lm;s}=e_{ls}Y_{lm;s}\text{ ,}
\end{equation*}%
where $\left\{ e_{ls}\right\} _{l=s,s+1}=\left\{ \left( l-s\right) \left(
l+s+1\right) \right\} _{l=s,s+1}$ is the associated sequence of eigenvalues
and $\left\{ Y_{lm;s}\right\} ,$ $l=s,s+1,...;m=-l,...,l$ is the sequence of
orthonormal spherical harmonics, which we define by
\begin{eqnarray*}
Y_{lm;s} &:&=\sqrt{\frac{(l-s)!}{(l+s)!}}\eth Y_{lm}\text{ for }s>0\text{ ,}
\\
Y_{lm;s} &:&=\sqrt{\frac{(l+s)!}{(l-s)!}}\eth Y_{lm}\text{ for }s<0\text{ .}
\end{eqnarray*}%
Again, as before it should be noted that in the spin case the operators
depend on the choice of the coordinates, differently from the scalar case.
As discussed by \cite{glm}, \cite{leosa}, \cite{mal} the spin construction
could be alternatively provided in terms of the so-called spin-weighted
representation of the special group of rotations $SO(3)$, indeed spin
spherical harmonics can be related to the so-called Wigner's matrices, see
again \cite{vmk},\cite{vilenkin}. In particular, it is then possible to show
that the spin spherical harmonics are themselves an orthonormal system, i.e.
they satisfy%
\begin{equation*}
\int_{\mathbb{S}^{2}}Y_{lm;s}\overline{Y}_{lm;s}dx=\int_{0}^{2\pi
}\int_{0}^{\pi }Y_{lm;s}(\vartheta ,\varphi )\overline{Y}_{lm;s}(\vartheta
,\varphi )\sin \vartheta d\vartheta d\varphi =\delta _{l}^{l^{\prime
}}\delta _{m}^{m^{\prime }}\text{ .}
\end{equation*}%
As for the scalar case,
\begin{equation*}
L_{s}^{2}\left( \mathbb{S}^{2}\right) =\bigoplus_{l=0}^{\infty }\mathcal{H}%
_{l}\qquad \mathcal{H}_{l}:=span\left\{ Y_{lm;s};m=-l,\ldots ,l\right\}
\text{ ,}
\end{equation*}%
and the following representation holds%
\begin{equation*}
F_{s}\left( x\right) =\sum_{l}\sum_{m}a_{lm;s}Y_{lm;s}(x)\text{{\ ,}}
\end{equation*}%
in the $L_{s}^{2}$ sense, i.e.
\begin{equation*}
\lim_{L\rightarrow \infty }\int_{\mathbb{S}^{2}}\left| F_{s}\left( x\right)
-\sum_{l=|s|}^{L}\sum_{m=-l}^{l}a_{lm;s}Y_{lm;s}(x)\right| ^{2}dx=0\text{ .}
\end{equation*}%
Here, the spherical harmonics coefficients $a_{lm;s}:=\int_{\mathbb{S}%
^{2}}F_{s}\overline{Y}_{lm}dx${\ are such that}%
\begin{equation*}
a_{lm;s}=a_{lm;E}+ia_{lm;M}\text{ ,}
\end{equation*}%
where $\left\{ a_{lm;E}\right\} ,\left\{ a_{lm;M}\right\} $ are the
coefficients of two standard (scalar-valued) spherical functions, which in
the physical literature are labelled the electric and magnetic components of
the spin function $F_{s},$ see again \cite{gelmar},\cite{gelmar2010} for
more discussion.

\section{Spin and Mixed Needlets}

\label{sec:spinneedlets}

\subsection{Definition}

We start by recalling the definition of scalar needlets, which were
introduced by \cite{npw1} and \cite{npw2} as:%
\begin{equation*}
\psi _{jk}\left( x\right) =\sqrt{\lambda _{jk}}\sum_{l}b\left( \frac{l}{B^{j}%
}\right) \sum_{m=-l}^{l}Y_{lm}\left( x\right) \overline{Y}_{lm}\left( \xi
_{jk}\right) \text{ , }\forall x\in \mathbb{S}^{2};
\end{equation*}%
here $\left\{ \xi _{jk},\lambda _{jk}\right\} $ are a set of
cubature points and weights ensuring that:
\begin{equation*}
\sum_{jk}\lambda _{jk}Y_{lm}\left( \xi _{jk}\right) \overline{Y}_{l^{\prime
}m^{\prime }}\left( \xi _{jk}\right) =\int_{\mathbb{S}^{2}}Y_{lm}\left(
x\right) \overline{Y}_{l^{\prime }m^{\prime }}\left( x\right) dx=\delta
_{l}^{l^{\prime }}\delta _{m}^{m^{\prime }},
\end{equation*}%
$b\left( \cdot \right) $ is a compactly supported $C^{\infty }$ function
satisfying the partition of unity property:
\begin{equation*}
\sum_{j}b^{2}(\frac{l}{B^{j}})\equiv 1
\end{equation*}%
for all $l\geq 1$ , and $B>1$ is a bandwidth parameter. For a fixed value of
$B$, we denote $\left\{ \mathcal{X}_{j}\right\} _{j\geq 0}$ the nested
sequence of cubature points corresponding to the space $\mathcal{K}_{\left[
2B^{j+1}\right] }$, where $\left[ \cdot \right] $ represents the integer
part and $\mathcal{K}_{L}=\oplus _{l=0}^{L}\mathcal{H}_{l}$ is the space
spanned by spherical harmonics up to order $L$. For each $j$, the cubature
points are almost distributed as an $\alpha _{j}$-net, with $\alpha
_{j}:=kB^{-j}$, the coefficients $\left\{ \lambda _{jk}\right\} $ are such
that $cB^{-2j}\leq \lambda _{jk}\leq CB^{-2j}$, with $c,C\in \mathbb{R}$,
and $N_{j}=card\left\{ \mathcal{X}_{j}\right\} \approx B^{2j}$, see for
instance \cite{bkmpBer} for more details.\newline
The construction of spin needlets (as provided by (\cite{gelmar})) is
formally similar to the scalar case, although as we discuss below it entails
deep differences in terms of the spaces involved. Indeed, spin needlets are
defined as follows:
\begin{equation}
\psi _{jk;s}\left( x\right) =\sqrt{\lambda _{jk}}\sum_{l}b\left( \frac{\sqrt{%
e_{l,s}}}{B^{j}}\right) \sum_{m=-l}^{l}\overline{Y}_{lm;s}\left( \xi
_{jk}\right) Y_{lm;s}\left( x\right) \text{ ,}  \label{nspin}
\end{equation}%
where $\left\{ \lambda _{jk},\xi _{jk}\right\} $ are, as before, cubature
weights and cubature points, $b\left( \cdot \right) \in C^{\infty }$ is
nonnegative, it is compactly supported in $\left[ 1/B,B\right] $ and
satisfies the partition of unity property. Note, however, that the
mathematical meaning of (\ref{nspin}) is rather different from the scalar
case; indeed $\psi _{jk;s}\left( x\right) $ is to be viewed as a spin $s$
function with respect to rotations of the tangent plane $\mathbb{T}_{x},$
and a spin $-s$ function with respect to rotations of the tangent plane $%
\mathbb{T}_{\xi _{jk}}.$ Moreover, as $Y_{lm;s}\left( \xi _{jk}\right)
,Y_{lm;s}\left( x\right) $ live on two different tangent planes $\mathbb{T}%
_{\xi _{jk}},\mathbb{T}_{x},$ the product $\overline{Y}_{lm;s}\left( \xi
_{jk}\right) Y_{lm;s}\left( x\right) $ is not defined and the notation $%
\overline{Y}_{lm;s}\left( \xi _{jk}\right) \otimes Y_{lm;s}\left( x\right) $
would be more appropriate. As a consequence, the spin needlet operators acts
on spin $s$ functions to produce spin $s$ coefficients
\begin{eqnarray}
\left\langle F_{s},\psi _{jk;s}\left( x\right) \right\rangle &=&\int_{%
\mathbb{S}^{2}}F_{s}(x)\overline{\psi }_{jk;s}\left( x\right) dx  \notag \\
&=&\sqrt{\lambda _{jk}}\sum_{lm}b\left( \frac{\sqrt{e_{l,s}}}{B^{j}}\right)
a_{lm;s}Y_{lm;s}\left( \xi _{jk}\right)  \notag \\
&=&:\beta _{jk;s}\text{ .}  \label{volley}
\end{eqnarray}%
Therefore, $\psi _{jk;s}$ induces the linear map (\ref{volley}) from spin $s$
quantities to spin $s$ wavelet coefficients $\beta _{jk;s}$, while in the
scalar case ($s=0$) needlets generate a linear map from scalar quantities to
scalar quantities. Indeed, if $u$ is a spin $s$ vector at $\xi _{jk}$, $\psi
_{jk;s}\left( x\right) u$ becomes a spin $s$ vector at $\xi _{jk}$, since
the product of spin $-s$ and spin $s$ vectors \textit{at a point x} is a
well-defined complex number, independently of the choice of coordinate
system.

To provide a clearer interpretation to the previous expression, recall the
decomposition of the functional space $L_{s}^{2}\left( \mathbb{S}^{2}\right)
=\bigoplus_{l\geq 0}\mathcal{H}_{l}$. We can hence define the following
operators on $\mathcal{H}_{l}$:
\begin{eqnarray*}
K_{j}\left( x,y\right) &=&\sum_{l}b^{2}\left( \frac{\sqrt{e_{l,s}}}{B^{j}}%
\right) Y_{lm;s}\left( x\right) \overline{Y}_{lm;s}\left( y\right) \\
\Lambda _{j}\left( x,y\right) &=&\sum_{l}b\left( \frac{\sqrt{e_{l,s}}}{B^{j}}%
\right) Y_{lm;s}\left( x\right) \overline{Y}_{lm;s}\left( y\right)
\end{eqnarray*}%
such that the reproducing kernel property holds:
\begin{equation*}
\int_{\mathbb{S}^{2}}\Lambda _{j}\left( x,y\right) \overline{\Lambda }%
_{j}\left( y,z\right) dy=K_{j}(x,z)\text{ .}
\end{equation*}%
Spin needlets can be derived by discretizing this operator by using the
reproducing kernel property. In fact $\Lambda _{j}$ is such that:
\begin{equation*}
z\rightarrow \Lambda _{j}\left( x,z\right) \in \mathcal{K}_{\left[ B^{2j+1}%
\right] }\text{ ,}
\end{equation*}%
and therefore:
\begin{equation*}
z\rightarrow \Lambda _{j}\left( x,z\right) \overline{\Lambda }_{j}\left(
z,y\right) \in \mathcal{K}_{\left[ B^{4j+2}\right] }\text{ .}
\end{equation*}%
After discretization, we obtain:
\begin{equation*}
K_{j}\left( x,y\right) =\sum_{\xi _{jk}\in \mathcal{K}_{\left[ B^{4j+2}%
\right] }}\lambda _{jk}\Lambda _{j}\left( x,\xi _{jk}\right) \overline{%
\Lambda }_{j}\left( \xi _{jk},y\right) \text{ ,}
\end{equation*}%
where we exploit the fact that the pairs $\left\{ \lambda _{jk},\xi
_{jk}\right\} $ can be chosen to form exact cubature points and
weights (\cite{bkmp09}). Then
\begin{eqnarray*}
K_{j}f\left( x\right) &=&\int_{\mathbb{S}^{2}}K_{j}\left( x,y\right) f\left(
y\right) dy \\
&=&\int_{\mathbb{S}^{2}}{\sum_{\xi _{jk}\in \mathcal{K}_{\left[B^{4j+2}%
\right] }}}\lambda _{jk}\Lambda _{j}\left( x,\xi _{jk}\right) \overline{%
\Lambda }_{j}\left( \xi _{jk},y\right) f\left( y\right) dy \\
&=&\sum_{\xi _{jk}\in \mathcal{K}_{\left[B^{4j+2}\right] }}\sqrt{\lambda
_{jk}}\Lambda _{j}\left( x,\xi _{jk}\right) \int_{\mathbb{S}^{2}}\sqrt{%
\lambda _{jk}}\overline{\Lambda }_{j}\left( \xi _{jk},y\right) f\left(
y\right) dy \\
&=&\sum_{\xi _{jk}\in \mathcal{K}_{\left[ B^{4j+2}\right] }}\beta
_{jk;s}\psi _{jk;s}\text{ ,}
\end{eqnarray*}%
where
\begin{equation*}
\psi _{jk;s}=\sqrt{\lambda _{jk}}\Lambda _{j}\left( x,\xi _{jk}\right) \text{
.}
\end{equation*}

As a minor point, note that for the argument of the function $b\left( \cdot
\right) $ we have used here the square root of $e_{l,s}$,the eigenvalue of
the corresponding spin spherical harmonics, while in the scalar case \cite%
{npw1},\cite{npw2} proposed to adopt $l$. However it is trivial to observe
that, for fixed $s$:
\begin{equation*}
\lim_{l\rightarrow \infty }\frac{\sqrt{e_{l,s}}}{l}=\lim_{l\rightarrow
\infty }\frac{\sqrt{\left( l-s\right) \left( l+s+1\right) }}{l}=1\text{ .}
\end{equation*}

\subsection{Some Properties}

We report some important properties for spin needlets, very similar to those
in scalar case (see \cite{npw1}, \cite{npw2}). Indeed, from the previous
discussion it follows easily that $\left| \psi _{jk;s}\right| ^{2}$ is a
well-defined scalar quantity. The following Localization property is hence
well-defined (see \cite{gelmar}): for any $M\in \mathbb{N}$, there exists a
constant $c_{M}>0$ such that for every $x\in \mathbb{S}^{2}$:%
\begin{equation*}
\left| \psi _{jk;s}\left( x\right) \right| \leq \frac{c_{M}B^{j}}{\left(
1+B^{j}\arccos \left( \langle \xi _{jk},x\rangle \right) \right) ^{M}}\text{
.}
\end{equation*}

Let us recall from (\ref{volley}) that%
\begin{equation*}
\beta _{jk;s}=\int_{\mathbb{S}^{2}}F_{s}\left( x\right) \overline{Y}%
_{jk;s}\left( x\right) dx=\sqrt{\lambda _{jk}}\sum_{l}b\left( \frac{\sqrt{%
e_{ls}}}{B^{j}}\right) \sum_{m=-l}^{l}a_{lm;s}Y_{lm;s}\left( \xi
_{jk}\right) \text{ , }
\end{equation*}%
and the following \textit{reconstruction formula} holds:
\begin{equation*}
F_{s}\left( x\right) =\sum_{j}\sum_{k}\beta _{jk;s}\psi _{jk;s}\left(
x\right) \text{ .}
\end{equation*}%
It is simple to check that the squared coefficients $\left| \beta
_{jk;s}\right| ^{2}$ following quantities are scalar. In the following, we
will need both the $L_{s}^{2}\left( \mathbb{S}^{2}\right) $ and the $%
L_{s}^{p}\left( \mathbb{S}^{2}\right) $ norm of $\psi _{jk;s}$. Let us start
by observing that:%
\begin{eqnarray*}
\left\| \psi _{jk;s}\right\| _{L_{s}^{2}\left( \mathbb{S}^{2}\right) }^{2}
&=&\lambda _{jk}\sum_{l}b^{2}\left( \frac{\sqrt{e_{ls}}}{B^{j}}\right)
\sum_{m=-l}^{l}Y_{lm;s}\left( \xi _{jk}\right) \overline{Y}_{lm;s}\left( \xi
_{jk}\right) \int_{\mathbb{S}^{2}}\ Y_{lm;s}\left( x\right) \overline{Y}%
_{lm;s}\left( x\right) dx= \\
&=&\lambda _{jk}\sum_{l=B^{j-1}}^{B^{j+1}}b^{2}\left( \frac{\sqrt{e_{ls}}}{%
B^{j}}\right) \sum_{m=-l}^{l}Y_{lm;s}\left( \xi _{jk}\right) \overline{Y}%
_{lm;s}\left( \xi _{jk}\right)  \\
&=&\lambda _{jk}\sum_{l=B^{j-1}}^{B^{j+1}}\frac{2l+1}{4\pi }b^{2}\left(
\frac{\sqrt{e_{ls}}}{B^{j}}\right) =:\tau _{jk;s}^{2}\text{ .}
\end{eqnarray*}%
As discussed by \cite{bkmpBer}, \cite{bkmp09}, \cite{gelmar}, there exist
positive constants $c_{1},c_{2}$ such that $c_{1}N_{j}^{-1}\leq \lambda
_{jk}\leq c_{2}N_{j}^{-1}$.

Throughout the rest of the paper, to simplify notations we shall assume to
be dealing with sections of line bundles such that $F_{s}=(I-P_{s})F_{s},$ $%
P_{s}$ denoting the projection operator on the $s$ spin spherical
harmonics. In other words, the component at $l=s$ is assumed to be
null; from the point of view of motivating applications, this is a
very reasonable assumption, indeed for polarization or weak lensing
experiments the so-called quadrupole term $l=s=2$ has no physical
meaning. The situation is indeed analogous to the standard scalar
case, where the constant term $s=0$ cannot even be measured by
ongoing (so-called \emph{differential}) experiments. Under these
circumstances, as shown in \cite{bkmp09}, spin needlets make up a
tight frame system, i.e. for all $F_{s}\in
L_{s}^{2}(\mathbb{S}^{2})$ ,
\begin{equation*}
\left\| F_{s}\right\| _{L_{s}^{2}\left( \mathbb{S}^{2}\right)
}^{2}=\sum_{jk}\left| \beta _{jk;s}\right| ^{2},
\end{equation*}%
whence we have easily
\begin{equation*}
\sum_{jk}\left| \langle \psi _{j_{1},k_{1};s},\psi _{jk:s}\rangle \right| {%
^{2}}=\left\| \psi _{j_{1},k_{1};s}\right\| _{L_{s}^{2}\left( \mathbb{S}%
^{2}\right) }^{4}+\sum_{j\neq j_{1},k\neq k_{1}}\left| \langle \psi
_{j_{1},k_{1};s}|\psi _{jk:s}\rangle \right| {^{2}}\leq \left\| \psi
_{j_{1},k_{1};s}\right\| _{L_{s}^{2}\left( \mathbb{S}^{2}\right) }^{2},
\end{equation*}%
whence%
\begin{equation*}
\left\| \psi _{jk;s}\right\| _{L_{s}^{2}\left( \mathbb{S}^{2}\right) }\leq 1%
\text{ .}
\end{equation*}%
More generally, it is shown in \cite{bkmp09}, \cite{gelmar2010} that for all
$1\leq p\leq \infty ,$ there exist positive constants $c_{p},C_{p}$ such that%
\begin{equation}
c_{p}B^{2j\left( \frac{1}{2}-\frac{1}{p}\right) }\leq \left\| \psi
_{jk;s}\right\| _{L_{s}^{p}\left( \mathbb{S}^{2}\right) }=\left( \int_{%
\mathbb{S}^{2}}\left| \psi _{jk;s}\right| ^{p}dx\right) ^{\frac{1}{p}}\leq
C_{p}B^{2j\left( \frac{1}{2}-\frac{1}{p}\right) }.  \label{normbound}
\end{equation}

\subsection{Mixed Needlets and their properties}

Mixed Needlets were introduced in \cite{gelmar2010}; they are defined as
\begin{equation*}
\psi _{jk;s\mathcal{M}}\left( x\right) =\sqrt{\lambda _{jk}}\sum_{l\geq
\left| s\right| }b\left( \frac{\sqrt{e_{ls}}}{B^{j}}\right)
\sum_{m}Y_{lm;s}\left( x\right) \overline{Y}_{lm}\left( \xi _{jk}\right)
\text{ ,}
\end{equation*}%
with corresponding needlet coefficients%
\begin{equation*}
\beta _{jk;s\mathcal{M}}=\int_{\mathbb{S}^{2}}\overline{\psi }_{jk;s\mathcal{%
M}}\left( x\right) F_{s}\left( x\right) dx\text{ .}
\end{equation*}%
Mixed needlets form a tight frame system, with the same set of cubature
points and weights as for the scalar case, $\left\{ \xi _{jk},\lambda
_{jk}\right\} .$ When $F_{s}\in L_{s}^{2}\left( \mathbb{S}^{2}\right) ,$ we
have also%
\begin{equation*}
\beta _{jk;s\mathcal{M}}=\sqrt{\lambda _{jk}}\sum_{l\geq \left| s\right|
}b\left( \frac{\sqrt{e_{ls}}}{B^{j}}\right) \sum_{m}a_{lm;s}Y_{lm}\left( \xi
_{jk}\right) \text{ ,}
\end{equation*}%
and mixed needlets form a tight frame system. It should be noted that the
coefficients $\left\{ \beta _{jk;s\mathcal{M}}\right\} $ are scalar,
complex-valued random variables, indeed for square integrable sections we
have
\begin{eqnarray*}
\beta _{jk;s\mathcal{M}} &=&\sqrt{\lambda _{jk}}\sum_{l\geq \left| s\right|
}b\left( \frac{\sqrt{e_{ls}}}{B^{j}}\right) \sum_{m}\left\{
a_{lm;E}+ia_{lm;M}\right\} Y_{lm}\left( \xi _{jk}\right)  \\
&=&:\beta _{jk;E}+i\beta _{jk;M}\text{ ,}
\end{eqnarray*}%
where $\beta _{jk;E},\beta _{jk;M}$ could be viewed as the scalar needlet
coefficients of standard square integrable functions on the sphere. For
general $F_{s}\in L_{s}^{p}\left( \mathbb{S}^{2}\right) $ the reconstruction
formula holds, in the $L_{s}^{p}$ sense:
\begin{equation*}
F_{s}=\sum_{j}\sum_{k}\beta _{jk;s\mathcal{M}}\psi _{jk;s\mathcal{M}}\text{ .%
}
\end{equation*}%
Other properties of mixed needlets are analogous to those for the pure spin
construction. In particular, note that scalar and pure spin needlets are
both constructed by a convolution of a smooth function $b(.)$ with
projection operators such as, for instance, $\sum_{m}Y_{lm}(x)\overline{Y}%
_{lm}(y)$, $\sum_{m}Y_{lm;s}(x)\overline{Y}_{lm;s}(y)$. \newline
On the other mixed needlets are built by convolving $b(.)$ with $%
\sum_{m}Y_{lm}(x)\overline{Y}_{lm;s}(y),$ which is not a projection operator
(indeed $\sum_{m}Y_{lm}(x)\overline{Y}_{lm;s}(x)\equiv 0).$ It comes
therefore to some extent as a surprise that mixed needlets do indeed enjoy
localization properties, indeed we have (see again \cite{gelmar2010}): for
each $M>0$ there exists a constant $C_{M}$ such that:
\begin{equation*}
\left| \psi _{jk;s\mathcal{M}}\right| \leq \frac{C_{M}B^{j}}{\left(
1+B^{j}\arccos \left( \langle x,\xi _{jk}\rangle \right) \right) ^{M}}\text{
.}
\end{equation*}%
Building upon this localization property, it is indeed possible to establish
the following bounds (see for more details \cite{gelmar2010}):
\begin{equation}
c_{1}B^{2j\left( \frac{1}{2}-\frac{1}{p}\right) }\leq \left\| \psi _{jk;s%
\mathcal{M}}\right\| _{L_{s}^{p}\left( \mathbb{S}^{2}\right) }\leq
c_{2}B^{2j\left( \frac{1}{2}-\frac{1}{p}\right) },\text{ }c_{1},c_{2}>0\text{
.}  \label{eqn:pnorm}
\end{equation}%
These constraints on the $L_{s}^{p}$ norms will have the greatest importance
for our results to follow. Also, for positive constants $c_{3},c_{4}$ and
arbitrary coefficients $\lambda _{k}$ we have%
\begin{equation}
c_{3}\sum_{k}\left| \lambda _{k}\right| ^{p}\left\| \psi _{jk;s\mathcal{M}%
}\right\| _{L_{s}^{p}\left( \mathbb{S}^{2}\right) }^{p}\leq \left\|
\sum_{k}\lambda _{k}\psi _{jk;s\mathcal{M}}\right\| _{L_{s}^{p}\left(
\mathbb{S}^{2}\right) }^{p}\leq c_{4}\sum_{k}\left| \lambda _{k}\right|
^{p}\left\| \psi _{jk;s\mathcal{M}}\right\| _{L_{s}^{p}\left( \mathbb{S}%
^{2}\right) }^{p}\text{ .}  \label{pitagora}
\end{equation}

\begin{remark}
While the mathematical construction and the properties that can be developed
on the mixed needlets are very similar to the spin case, there is a very
relevant difference among these approaches that will be very important for
our purposes. While, as we have already seen, $\psi _{jk;s}$ is formed by a
tensorial product among two terms belonging to two different spaces of spin $%
-s$ and $s$ such that $\beta _{jk_{s}}$ belongs to the spin $s$ space, $\psi
_{jk;s\mathcal{M}}$ induces a linear map from a spin $s$ vector at $\xi
_{jk} $ to a scalar (spin $0$) quantity, such that for a spin $s$ quantity $%
u $, the product $\overline{\psi }_{jk;s\mathcal{M}}\cdot u$ is always a
scalar quantity. .
\end{remark}

\section{Spin Besov spaces}

Our aim in this Section is to recall the definition of spin Besov spaces in
terms of approximation properties. These definitions and their
characterizations were provided by \cite{bkmp09}, \cite{gelmar2010}, to
which we refer for further details and discussion. Define first,
\begin{equation*}
G_{k}\left( F_{s},\pi \right) =\inf_{H_{s}\in \mathcal{H}_{k}}\left\|
F_{s}-H_{s}\right\| _{L_{s}^{\pi }\left( \mathbb{S}^{2}\right) }\text{ ,}
\end{equation*}%
i.e. the approximation error when replacing $F_s$ by an element in $\mathcal{%
H}_{k;s}$ . Then the Besov spin space $\mathcal{B}_{pq;s}^{r}$ is defined as
the space of functions such that $F_s\in L_{s}^{p}\left( \mathbb{S}%
^{2}\right) $ and%
\begin{equation*}
\left( \sum_{k=0}^{\infty }\frac{1}{k}\left( k^{r}G_{k}\left( F_s,\pi
\right) \right) ^{q}\right) <\infty \text{ .}
\end{equation*}%
As usual, the last condition can be easily shown to be equivalent to%
\begin{equation*}
\left( \sum_{j=0}^{\infty }\left( B^{jr}G_{B^{j}}\left( F_s,\pi \right)
\right) ^{q}\right) <\infty \text{ .}
\end{equation*}%
Moreover, $F_{s}\in \mathcal{B}_{\pi q;s}^{r}$ if and only if, for every $%
j=1,2,\ldots $%
\begin{equation*}
\left( \sum_{k}\left( \left| \beta _{jk;s}\right| \left\| \psi
_{jk;s}\right\| _{L_{s}^{\pi }\left( \mathbb{S}^{2}\right) }\right) ^{\pi
}\right) ^{\frac{1}{\pi }}=\varepsilon _{j}B^{-jr}
\end{equation*}%
where $\varepsilon _{j}\in \ell _{q}$ and $B>1$. By defining the Besov norm
as follows,%
\begin{equation*}
\left\| F_{s}\right\| _{\mathcal{B}_{\pi q;s}^{r}}=\left\{
\begin{matrix}
\left\| F_{s}\right\| _{L_s^\pi\left(\mathbb{S}^2\right)}+\left[
\sum_{j}B^{jq\left( r+\frac{1}{2}-\frac{1}{\pi}\right) }\left\{
\sum_{k}\left| \beta _{jk;s}\right| ^{\pi }\right\} ^{\frac{q}{\pi }}\right]
^{\frac{1}{q}} & \mbox{ if  }q<\infty \\
\left\| F_{s}\right\| _{L_s^\pi\left(\mathbb{S}^2\right) }+\underset{j}{\sup
}B^{j\left( r+\frac{1}{2}-\frac{1}{\pi }\right) }\left\| \left( \beta
_{jk;s}\right) _{k}\right\| {_{\ell _{\pi}}} & \mbox{ if  }q=\infty%
\end{matrix}%
\right. ,
\end{equation*}%
we obtain that, if $\max \left( 0,1/\pi-1/q\right) <r$ and $\pi,q>1$, then%
\begin{equation*}
F_{s}\in \mathcal{B}_{\pi q;s}^{r}\Leftrightarrow \left\| F_{s}\right\| _{%
\mathcal{B}_{\pi q;s}^{r}}<\infty \text{ .}
\end{equation*}%
Besov spaces are characterized by come convenient embeddings, which (as
always in this literature) will play a crucial role in our proofs to follow.
More precisely, we have that, for $\pi _{1}\leq \pi _{2},$ $q_{1}\leq q_{2}$%
\begin{equation}
\mathcal{B}_{\pi q_{1};s}^{r}\subset \mathcal{B}_{\pi q_{2};s}^{r}\text{ },%
\text{ }\mathcal{B}_{\pi _{2}q;s}^{r}\subset \mathcal{B}_{\pi _{1}q;s}^{r}%
\text{ , }\mathcal{B}_{\pi _{1}q;s}^{r}\subset \mathcal{B}_{\pi _{2}q;s}^{r-%
\frac{1}{\pi _{1}}+\frac{1}{\pi _{2}}}.  \label{embeddings}
\end{equation}%
The proof of (\ref{embeddings}) is exactly the same as for the scalar case,
see \cite{bkmpAoSb}. In particular%
\begin{eqnarray*}
\mathcal{B}_{\pi _{1}q;s}^{r} &\subset &\mathcal{B}_{\infty \infty ;s}^{r-%
\frac{1}{\pi _{1}}}\Longrightarrow \sup_{k}\left| \beta _{jk;s}\right|
\left\| \psi _{jk}\right\| _{L_s^\infty\left(\mathbb{S}^2\right)}=%
\varepsilon _{j}B^{-j(r-\frac{1}{\pi _{1}})} \\
&\Rightarrow &B^{j(r+1-\frac{1}{\pi _{1}})}\sup_{k}\left| \beta
_{jk;s}\right| <\infty \\
&\Rightarrow &B^{j}\sup_{k}\left| \beta _{jk;s}\right| <\infty \text{ .}
\end{eqnarray*}

\section{Nonparametric Regression on Spin Fiber Bundles}

\label{sec:estimators}

\subsection{The Regression Model}

We start by recalling the regression formula (\ref{eqn:regression}):
\begin{equation*}
Y_{i;s}=F_{s}\left( X_{i}\right) +\varepsilon _{i;s}\text{ .}
\end{equation*}%
Throughout this paper, we shall also assume that $\sup_{x}\left|
F_{s}(x)\right| =M<\infty $ . As discussed in the Introduction, we envisage
a situation where it is possible to collect data which can be viewed as
measurements on a spin fiber bundles, i.e. for instance the polarization of
the Cosmic Microwave Background (see \cite{kks96},\cite{selzad}, \cite{ck05}%
, \cite{glm}, \cite{ghmkp}), or the Weak Gravitational Lensing
effect on the images of distant Galaxies (see \cite{bridle}). To fix
ideas, we focus on this second example. As discussed for instance in
\cite{dode2004}, the gravitational \emph{shear }effect may be
loosely described as gravity transforming into a more elliptical
shape the image of galaxies. Of course the measurement of this shear
is subject to an experimental error, for instance because of the
unknown intrinsic ellipticity of the observed galaxy. Likewise, the
weak gravitational lensing may produce an alignment in the
inclination of nearby observations, but again this could be brought
in by random fluctuations. We refer to \cite{bridle,} for much more
detailed discussion on motivations and related challenges, which
currently involve huge amount of physicists; major satellite
experiments are at the planning stage, such as Euclid, see for
instance http://hetdex.org/other\_projects/euclid.php. To model the
above discussed framework, we introduce random directions of
observations $\left\{ X_{i}\in \mathbb{S}^{2}\right\} ,$ which we
take to be uniformly sampled over the
sky, and observational errors $\left\{ \varepsilon _{i;s}\right\} $, $%
i=1,2,...,n$; the latter are independent and identically distributed spin $s$
random variables, which we assume to be invariant in law with respect to
rotations in the tangent plane:
\begin{equation}
\varepsilon _{i;s}^{\prime }\overset{d}{=}\varepsilon _{i;s}e^{is\psi },%
\text{ for all }\psi \in \left[ 0,2\pi \right] \text{ , }i=1,2,...,n,
\label{manni}
\end{equation}%
$\overset{d}{=}$ denoting equality in law. As in \cite{bm}, (\ref{manni})
implies that
\begin{equation}
\func{Re}\varepsilon _{i;s}\overset{d}{=}\func{Im}\varepsilon _{i;s}\overset{%
d}{=}\widetilde{\varepsilon }_{i}\text{ .}  \label{manni2}
\end{equation}%
From (\ref{manni}), (\ref{manni2}) we have immediately%
\begin{equation*}
E\left[ \varepsilon _{i;s}\right] =E\left[ \func{Re}\varepsilon _{i;s}+i%
\func{Im}\varepsilon _{i;s}\right] =0\text{ ,}
\end{equation*}%
\begin{equation*}
Var\left( \varepsilon _{i;s}\right) =E\left| \varepsilon _{i;s}\right|
^{2}=2E\widetilde{\varepsilon }_{i}^{2}=:\sigma _{\varepsilon }^{2}\text{ .}
\end{equation*}%
Moreover, we shall assume that $\left\{ \widetilde{\varepsilon }_{i}\right\}
$ follows a sub-Gaussian distribution (ref.\cite{buko}), i.e. there exists a
number $a\geq 0$ such that for all $\lambda \in \mathbb{R}$ the following
inequality holds:
\begin{equation}
E\left[ e^{\lambda \widetilde{\varepsilon }_{i}}\right] \leq e^{\left( \frac{%
a^{2}\lambda ^{2}}{2}\right) }\text{ .}  \label{manni3}
\end{equation}%
We also define the \textit{sub-Gaussian standard }of the random variable $%
\widetilde{\varepsilon }_{i}$ as:
\begin{equation*}
\tau \left( \widetilde{\varepsilon }_{i}\right) = \inf \left\{ a\geq 0:E%
\left[ e^{\lambda \widetilde{\varepsilon }_{i}}\right] \leq e^{\left( \frac{%
a^{2}\lambda ^{2}}{2}\right) },\lambda \in \mathbb{R}\right\} <\infty \text{
.}
\end{equation*}%
It is immediate to check ( see \cite{buko}) that:
\begin{equation*}
\tau \left( \widetilde{\varepsilon }_{i}\right) =\sup_{\lambda \neq 0}\left[
\frac{2\log \left( E\left[ e^{\lambda \widetilde{\varepsilon }_{i}}\right]
\right) }{\lambda ^{2}}\right] ^{\frac{1}{2}},\text{ }E\left[ e^{\lambda X}%
\right] \leq e^{\frac{\lambda ^{2}\tau ^{2}\left( \widetilde{\varepsilon }%
_{i}\right) }{2}}.
\end{equation*}%
As well-known, a random variable is sub-Gaussian if and only if the moment
generating function is majorized by the moment generating function of a
zero-mean Gaussian random variable, whence the name \textit{sub-Gaussian}.
Indeed, the class of sub-Gaussian random variables contains, apart from the
Gaussian themselves, all bounded zero-mean random variables and, more
generally, all those random variables whose distribution tails decrease no
slower than the tails of the Gaussian. We recall the following, simple
results, whose proofs are available in \cite{buko}:

\begin{lemma}
\label{subG} \textit{Moment characterization for subGaussian random
variables.}\newline
Let $\tilde{\varepsilon}$ be a subGaussian random variable such that $%
E\left( \tilde{\varepsilon}\right) =0.$ We have that $E\left( \left( \tilde{%
\varepsilon}\right) ^{2}\right) \leq \tau \left( \tilde{\varepsilon}\right) $
and for all $p>0$ $E\left( \left| \tilde{\varepsilon}\right| ^{p}\right)
<\infty $.

%\item \textit{Hoeffding's inequality for subGaussian random variables.}%
%\newline
%Assume that $\tilde{\varepsilon}_{i}$, $i=1,\ldots ,n$ is a sequence of
%subGaussian random variables with corresponding standard $\tau _{i}$. The
%following inequality holds:%
%\begin{equation*}
%\mathbb{P}\left( \left| \sum_{i=1}^{n}\tilde{\varepsilon}_{i}\right|
%>x\right) \leq 2\exp \left( -\frac{x^{2}}{2\sum_{i=1}^{n}\tau \left(
%\widetilde{\varepsilon }_{i}\right) }\right) \text{ .}
%\end{equation*}
%\end{enumerate}
\end{lemma}

In view of Lemma (\ref{subG}), subGaussian random variables enjoy the same
moment inequalities and concentration properties as Gaussian or bounded
ones, and hence allow the implementation of the main technical tools in the
proofs of our asymptotic results to follow. In this sense, they seem to
provide a natural general framework for the analysis we must pursue.

\subsection{The estimation procedure}

The procedure we are going to investigate can be viewed as a form of needlet
thresholding in the spin fiber bundles case (we refer to \cite{bkmpAoSb} for
a similar approach, in the case of density estimation for standard scalar
directional data). As discussed in the Introduction, we have now two
alternative forms of needlets construction for the spin case, i.e. the pure
spin needlets of \cite{gelmar} and the mixed spin needlets of \cite%
{gelmar2010}. Our approach could be implemented for both techniques, and
indeed the proofs would be nearly identical. For definiteness, we shall
focus on the mixed needlets constructions, which yields coefficients which
are standard, complex-valued variables. For brevity's sake, however, we drop
the subscript $\mathcal{M}$. We start by defining, as usual, an unbiased
estimator for needlet coefficients. More precisely, we define
\begin{equation*}
\widehat{\beta }_{jk;s}:=\frac{1}{n}\sum_{i=1}^{n}Y_{i}\overline{\psi }%
_{jk;s}\left( X_{i}\right) \text{ , }i=1,2,...,n\text{ .}
\end{equation*}%
We have immediately:
\begin{eqnarray}
E\left( \widehat{\beta }_{jk;s}\right) &=&\frac{1}{n}\sum_{i=1}^{n}E\left[
\overline{\psi }_{jk;s}\left( X_{i}\right) F_{s}\left( X_{i}\right) +%
\overline{\psi }_{jk;s}\left( X_{i}\right) \varepsilon _{i;s}\right] =
\notag \\
&=&\int_{\mathbb{S}^{2}}\overline{\psi }_{jk;s}\left( X_{i}\right)
F_{s}\left( X_{i}\right) =\beta _{jk;s}\text{ .}  \label{meanbeta}
\end{eqnarray}%
Moreover%
\begin{eqnarray}
Var\left( \widehat{\beta }_{jk;s}\right) &=&Var\left( \frac{1}{n}%
\sum_{i=1}^{n}\overline{\psi }_{jk;s}\left( X_{i}\right) F_{s}(X_{i})+\frac{1%
}{n}\sum_{i=1}^{n}\overline{\psi }_{jk;s}\left( X_{i}\right) \varepsilon
_{i;s}\right) =  \label{varbeta} \\
&=&\frac{1}{n^{2}}\sum_{i=1}^{n}Var\left( \overline{\psi }_{jk;s}\left(
X_{i}\right) F_{s}(X_{i})\right) +\frac{1}{n^{2}}\sum_{i=1}^{n}{Var\left(
\overline{\psi }_{jk;s}\left( X_{i}\right) \varepsilon _{i;s}\right) }\text{
.}  \notag
\end{eqnarray}%
Now
\begin{equation*}
\frac{1}{n^{2}}\sum_{i=1}^{n}Var\left( \overline{\psi }_{jk;s}\left(
X_{i}\right) \varepsilon _{i;s}\right) =\frac{1}{n}\sigma _{\varepsilon
}^{2}\left\| \psi _{jk;s}\right\| _{L_{s}^{2}\left( \mathbb{S}^{2}\right)
}^{2}=\frac{1}{n}\sigma _{\varepsilon }^{2}\tau _{j}^{2}=:\frac{1}{n}\sigma
_{1\varepsilon ,j}^{2}
\end{equation*}%
where in the last equality we used the independence of the $\varepsilon
_{i;s}$. Note that obviously $\sigma _{1\varepsilon ,j}^{2}\leq \sigma
_{\varepsilon }^{2}$. Also%
\begin{equation*}
0\leq \frac{1}{n^{2}}\sum_{i=1}^{n}Var\left( \overline{\psi }_{jk;s}\left(
X_{i}\right) F_{s}(X_{i})\right) =\frac{1}{n}\int_{\mathbb{S}^{2}}\left|
\overline{\psi }_{jk;s}\left( x\right) F_{s}(x)\right| ^{2}dx\leq \frac{M^{2}%
}{n}\text{ ,}
\end{equation*}%
and we define $\sigma _{\varepsilon ,j}^{2}:=\sigma _{1\varepsilon ,j}^{2}+%
\frac{M^{2}}{n}$ . We then proceed with the (now classical) hard
thresholding procedure (see for instance \cite{donoho1}, \cite{WASA} and %
\cite{cai}). In particular, we fix the threshold as%
\begin{equation}
\kappa t_{n}=\kappa \sqrt{\frac{\log n}{n}}\text{ , }  \label{tn}
\end{equation}%
where $\kappa $ is a real positive constant, whose value will be discussed
later. Hence we define as usual%
\begin{equation}
\beta _{jk;s}^{\ast }=w_{jk}\widehat{\beta }_{jk;s}\text{ , }w_{jk}=\mathbb{I%
}_{\left\{ |\widehat{\beta }_{jk;s}|>\kappa t_{n}\right\} }\text{ , }
\label{thresho}
\end{equation}%
where $\mathbb{I}_{A}$ denotes as usual the indicator function of the set $A$%
. The thresholding estimator is hence%
\begin{equation}
F_{s}^{\ast }\left( x\right) =\sum_{j=1}^{J_{n}}\sum_{k=1}^{N_{j}}\beta
_{jk;s}^{\ast }\psi _{jk;s}\left( x\right) \text{ . }  \label{basest}
\end{equation}%
In \ref{basest}, $\ J_{n}$ represents a cut-off frequency, which we shall
fix at $B^{J_{n}}=\sqrt{\frac{n}{\log n}},$ whereas $N_{j}$ is the
cardinality of the cubature point set at frequency $j;$ it is known (see for
instance \cite{bkmpBer}) that there exist positive constants $c_{1},c_{2}$
such that $c_{1}B^{2j}\leq N_j\leq c_{2}B^{2j}$ (written $N^{j}\approx
B^{2j})$. Our main result is to show that thresholding estimates achieve
'nearly optimal' (up to logarithmic factors) rates with respect to general $%
L_{s}^{p}\left( \mathbb{S}^{2}\right) $ loss functions..

\begin{theorem}
\label{th:loss} Let $F_{s}\in \mathcal{B}_{\pi q;s}^{r}(G),$ the Besov ball
such that $\left\| F_{s}\right\| _{\mathcal{B}_{\pi q;s}^{r}}\leq G<\infty ,$
$r-\frac{2}{\pi }>0,$ and consider $F_{s}^{\ast }$ defined by (\ref{basest}, %
\ref{tn}, \ref{thresho}). For $1\leq p<\infty ,$ there exist $\kappa >0$
such that we have%
\begin{equation*}
\sup_{F_{s}\in \mathcal{B}_{\pi q;s}^{r}}E\left\| F_{s}^{\ast
}-F_{s}\right\| _{L_{s}^{p}}^{p}\leq C_{p}\left\{ \log n\right\} ^{p}\left[
\frac{n}{\log n}\right] ^{-\alpha (r,\pi ,p)},
\end{equation*}%
\begin{equation*}
\alpha (r,\pi ,p)=\left\{
\begin{array}{c}
\frac{rp}{2r+2}\text{ for }\pi \geq \frac{2p}{2r+2}\text{ (regular zone)} \\
\frac{p(r-2(\frac{1}{\pi }-\frac{1}{p}))}{2(r-2(\frac{1}{\pi }-\frac{1}{2}))}%
\text{ for }\pi \leq \frac{2p}{2r+2}\text{ (sparse zone)}%
\end{array}%
\right. .
\end{equation*}%
Also, for $p=\infty $%
\begin{equation*}
\sup_{F_{s}\in \mathcal{B}_{\pi q;s}^{r}}E\left\| F_{s}^{\ast
}-F_{s}\right\| _{L_{s}^{\infty }}\leq C_{\infty }\left[ \frac{n}{\log n}%
\right] ^{-\alpha (r,\pi ,\infty )},\text{ }\alpha (r,\pi ,\infty )=\frac{(r-%
\frac{2}{\pi })}{2(r-2(\frac{1}{\pi }-\frac{1}{2}))}\text{ .}
\end{equation*}
\end{theorem}

\begin{remark}
The definitions of regular and sparse zones are classical, and so are the
rates we obtained, which indeed correspond (for instance) to those presented
by \cite{bkmpAoSb}. For brevity's sake, we do not prove that these rates are
indeed minimax (up to logarithmic terms), but it seems easy to achieve this
goal by application of classical arguments, as for instance presented by %
\cite{thres}. It is trivial to note that for $\pi =\frac{2p}{2r+2}=\frac{p}{%
r+1}$ we have%
\begin{eqnarray*}
\frac{(r-2(\frac{1}{\pi }-\frac{1}{p}))}{2(r-2(\frac{1}{\pi }-\frac{1}{2}))}
&=&\frac{(r-2(\frac{r+1}{p}-\frac{1}{p}))}{2(r-2(\frac{r+1}{p}-\frac{1}{2}))}%
=\frac{r(p-2)}{2r(p-2)+2(p-2)} \\
&=&\frac{r}{2r+2}\text{ ,}
\end{eqnarray*}%
and also
\begin{eqnarray*}
\frac{rp}{2r+2} &\geq &\frac{p(r-2(\frac{1}{\pi }-\frac{1}{p}))}{2(r-2(\frac{%
1}{\pi }-\frac{1}{2}))}\text{ in the regular zone ,} \\
\frac{rp}{2r+2} &\leq &\frac{p(r-2(\frac{1}{\pi }-\frac{1}{p}))}{2(r-2(\frac{%
1}{\pi }-\frac{1}{2}))}\text{ in the sparse zone .}
\end{eqnarray*}%
Of course $\alpha (r,\pi ,p)<\frac{1}{2}$ , $\lim_{r\rightarrow \infty
}\alpha (r,\pi ,p)=\frac{1}{2}$ .
\end{remark}

\begin{remark}
For $s=0,$ our results cover adaptive nonparametric regression for
complex-valued, scalar functions. Again, the rates correspond to the usual
nearly minimax bounds.
\end{remark}

The proof of Theorem \ref{th:loss} is provided in the Section to follow.

\section{Proofs}

Our arguments will follow classical approaches in this area, as presented
for instance by \cite{thres}, \cite{bkmpAoSb}.

\subsection{An auxiliary result}

We shall need in the sequel some sharp bounds which are provided in the
following result. The arguments are close, for instance, to those for the
inequality (65) on page 1088 of \cite{thres} where the case of a scalar
Gaussian noise is considered: see also Proposition 15 in \cite{bkmpAoSb}.

\begin{proposition}
\label{prop:1} Let $\left\{ \varepsilon _{i;s}\right\} $ be such that (\ref%
{manni}) and (\ref{manni3}) are fulfilled. Assume also that $M:=\left\|
F_{s}\right\| _{\infty }<\infty $. For all $\gamma >0$ and for all $j$ such
that $B^{j}\leq \sqrt{n/\log n}$, there exists $\kappa _{\gamma }>0$ such
that for $\kappa >\kappa _{\gamma }$ the following inequality holds:
\begin{equation}
\mathbb{P}\left( \frac{1}{n}\left| \sum_{i=1}^{n}{\ \overline{\psi }%
_{jk;s}\varepsilon _{i;s}}\right| >\kappa \sqrt{\frac{\log n}{n}}\right)
\leq Cn^{-\gamma }.  \label{eqn:prop1}
\end{equation}%
where $\gamma \approx \kappa ^{4/3}$. Moreover, for all $p>0$ we have
\begin{eqnarray}
&&E\left[ \left| \widehat{\beta }_{jk;s}-\beta _{jk;s}\right| ^{p}\right]
\leq C_{p}n^{-\frac{p}{2}}  \label{eqn:prop2} \\
&&E\left[ \sup_{k}\left| \widehat{\beta }_{jk;s}-\beta _{jk;s}\right| ^{p}%
\right] \leq C_{\infty }(j+1)^{p-1}n^{-p/2}.  \label{eqn:prop3}
\end{eqnarray}
\end{proposition}

\begin{remark}
It is possible to obtain sharp analytic expressions for $\kappa
,C_{p},C_{\infty },$ for instance by arguing as in Lemma 16 of \cite%
{bkmpAoSb}.
\end{remark}

\textbf{Proof} Note first that
\begin{eqnarray}
&&\widehat{\beta }_{jk;s}=\frac{1}{n}\sum_{i=1}^{n}\overline{\psi }%
_{jk;s}\left( F_{s}\left( X_{i}\right) +\varepsilon _{i;s}\right)
\label{eqn:split} \\
&&\beta _{jk;s}=E\left( \widehat{\beta }_{jk;s}\right) =\frac{1}{n}%
\sum_{i=1}^{n}{\ E\left( \overline{\psi }_{jk;s}\left( X_{i}\right)
F_{s}\left( X_{i}\right) \right) }\text{ ,}
\end{eqnarray}%
and%
\begin{eqnarray*}
\widehat{\beta }_{jk;s}-\beta _{jk;s} &=& \frac{1}{n}\sum_{i=1}^{n}\left\{
\left( \overline{\psi }_{jk;s}F_{s}\left( X_{i}\right) \right) {-}E\left(
\overline{\psi }_{jk;s}F_{s}\left( X_{i}\right) \right) \right\} +\frac{1}{n}%
\sum_{i=1}^{n}{\overline{\psi }}_{jk;s}\varepsilon _{i;s} \\
&=& \frac{1}{n}\sum_{i=1}^n \Psi _{jk;s}(X_{i}) + \frac{1}{n}\sum_{i=1}^n{%
\overline{\psi }}_{jk;s}\varepsilon _{i;s}\text{ .}
\end{eqnarray*}%
where
\begin{equation*}
\Psi _{jk;s}(X_{i}):=\overline{\psi }_{jk;s}(X_{i})F_{s}\left( X_{i}\right)
-E\left( \overline{\psi }_{jk;s}(X_{i})F_{s}\left( X_{i}\right) \right)\text{
.}
\end{equation*}%
%
%
%We can split the first term in (\ref{eqn:prop1}) into two addenda:%
Consider $\mathbb{P}_{\beta }\left( x\right) :=\mathbb{P}\left( \left|
\widehat{\beta }_{jk;s}-\beta _{jk;s}\right| >x\right) $:
\begin{equation}
\mathbb{P}_{\beta }\left( x\right) \leq \mathbb{P}_{F}\left( x\right) +%
\mathbb{P}_{\varepsilon }\left( x\right)  \label{pe1}
\end{equation}%
where:
\begin{eqnarray*}
&&\mathbb{P}_{F}\left( x\right) =\mathbb{P}\left( \frac{1}{n}\left|
\sum_{i=1}^{n}\Psi _{jk;s}(X_{i})\right| >\frac{1}{2}x\right) \text{ ,} \\
&&\mathbb{P}_{\varepsilon }\left( x\right) =\mathbb{P}\left( \frac{1}{n}%
\left| \sum_{i=1}^{n}\overline{\psi }_{jk;s}\varepsilon _{i;s}\right| >\frac{%
1}{2}x\right) \text{ .}
\end{eqnarray*}%
As before, we can split these sums into a real and imaginary part, to which
we can apply separately the following procedures for both real and imaginary
part in $\mathbb{P}_{F}\left( x\right)$ and $\mathbb{P}_{\varepsilon }\left(
x\right)$, that give the same results.\newline
As far as $\mathbb{P}_{F}\left( x\right)$ is concerned, we use the fact that
$\Psi _{jk;s}(X_{i}) $ are i.i.d random variables such that for each of
them:
\begin{eqnarray*}
&& \sup \left| \Psi _{jk;s}(X_{i}) \right| \leq 2cMB^{j} \\
&& E\left( \left| \Psi_{jk;s}(X_{i}) \right| ^{2}\right) \leq E\left( \left|
\overline{\psi }_{jk;s}\left( X_{i}\right) F_{s}\left(X_{i}\right) \right|
^{2}\right) \leq M^{2}\left\| \psi _{jk;s}\right\|_{L_{s}^{2}\left( \mathbb{S%
}^{2}\right) }^{2}\leq M^{2} \text{ .}
\end{eqnarray*}

We therefore apply Bernstein inequality: For a sequence of i.i.d. random
variables $\left\{ X_{i}\right\} _{i=1}^{n}$ such that $E\left[ X_{i}\right]
=0$, $\left| X_{i}\right| \leq M$ and $E\left[ X_{i}^{2}\right] =\sigma ^{2}$%
, we have%
\begin{equation}  \label{eqn:bern}
\mathbb{P}\left( \frac{1}{n}\left| \sum_{i=1}^{n}X_{i}\right| >x\right) \leq
2\exp \left( -\frac{nx^{2}}{2\left( \sigma ^{2}+\frac{1}{3}Mx\right) }%
\right) \text{ ,}
\end{equation}%
see for instance \cite{WASA}. for a proof.

By applying Bernstein, we obtain:
\begin{equation}  \label{eqn:PFx}
\mathbb{P}_{F}\left( x\right) \leq 4\exp \left( -\frac{n\frac{x^{2}}{4}}{%
2\left( M^{2}+\frac{1}{3}cMB^{j}x\right) }\right) \text{ ,}
\end{equation}%
where the value $4$ takes on count both real and imaginary parts.

%\begin{eqnarray*}
%\mathbb{P}\left( \left| \widehat{\beta }_{jk;s}-\beta _{jk;s}\right| >\kappa
%t_{n}\right) &\leq &\mathbb{P}\left( \frac{1}{n}\left| \sum_{i=1}^{n}%
%\overline{\psi }_{jk;s}F_{s}\left( X_{i}\right) -E\left( \overline{\psi }%
%_{jk;s}F_{s}\left( X_{i}\right) \right) \right| >\frac{\kappa }{2}\sqrt{%
%\frac{\log n}{n}}\right) + \\
%&+&\mathbb{P}\left( \frac{1}{n}\left| \sum_{i=1}^{n}\overline{\psi }%
%_{jk;s}\varepsilon _{i;s}\right| >\frac{\kappa }{2}\sqrt{\frac{\log n}{n}}%
%\right) =\mathbb{P}_{F}+\mathbb{P}_{\varepsilon }\text{ .}
%\end{eqnarray*}%

Fixing $x=\kappa t_n$, the following result is obtained:
\begin{equation*}
\mathbb{P}_{F}\left(\kappa t_n\right)\leq 4\exp \left( -\frac{n\left( \left(
k/2\right) \sqrt{\log n/n}\right) ^{2}}{\frac{2}{3}\left( 3M^{2}+cMB^{j}k%
\sqrt{\frac{\log n}{n}}\right) }\right) \text{ ,}
\end{equation*}%
and by choosing $j$ such that $B^{j}\leq \sqrt{\frac{n}{\log n}}$%
\begin{equation}  \label{eqn:PFtn}
\mathbb{P}_{F}\left(\kappa t_n\right)\leq 4\exp \left( -\frac{3k^{2}\log n}{%
8M\left( 3M+ck\right) }\right) =2n^{-\frac{3k^{2}}{8M\left( 3M+ck\right) }}.
\end{equation}

As far as $\mathbb{P}_{\varepsilon }\left(x\right)$ is concerned, consider
that conditionally on $\left( X_{1}^{\prime },\ldots ,X_{n}^{\prime }\right)
$, $\frac{1}{n}\sum_{i}{\overline{\psi }}_{jk;s}{\left( X_{i}^{\prime
}\right) }\varepsilon {_{i;s}} $ is a complex-valued subGaussian variable
with mean $0$ and variance $\frac{1}{n^{2}}\sum_{i=1}^{n}\left| \psi
_{jk;s}\left( X_{i}\right) \right| ^{2}\sigma _{\varepsilon }^{2}$.
Therefore, by using the Markov's inequality, we obtain:
\begin{eqnarray*}
\mathbb{P}_{\varepsilon }\left( x\right) &\leq & E\left( \left. \exp \left( -%
\frac{-nx^{2}}{\sigma _{\varepsilon} ^{2}\frac{8}{n}\sum_{i=1}^{n}{\left|
\psi _{jk;s}\right| ^{2}}}\right) \right| X_{1}^{\prime },\ldots
,X_{n}^{\prime }\right)
\end{eqnarray*}
Observe that $\left| \psi _{jk;s}\left( X_{i}^{\prime }\right) \right| ^{2}$
are i.i.d. variables bounded by $CB^{2j}$, such that $E\left( \left| \psi
_{jk;s}\left( X_{i}^{\prime }\right) \right| ^{2}\right) =\int_{\mathbb{S}%
^{2}}\left| \psi _{jk;s}\right| ^{2}dx=\left\| \psi _{jk;s}\right\|
_{L_{s}^{2}\left( \mathbb{S}^{2}\right) }^{2}\leq 1$. Therefore we split the
denominator into $2$ terms, using
\begin{equation*}
\mathbb{I}_{\left\{ \left| \frac{1}{n}\sum_{i=1}^{n}\left| \psi
_{jk;s}\left( X_{i}^{\prime }\right) \right| ^{2}{-}\left\| \psi
_{jk;s}\right\| _{L_{s}^{2}\left( \mathbb{S}^{2}\right) }^{2}\right| <\alpha
\right\} }\text{ and }\mathbb{I}_{\left\{ \left| \frac{1}{n}%
\sum_{i=1}^{n}\left| \psi _{jk;s}\left( X_{i}^{\prime }\right) \right| {^{2}-%
}\left\| \psi _{jk;s}\right\| _{L_{s}^{2}\left( \mathbb{S}^{2}\right)
}^{2}\right| \geq \alpha \right\} }\text{ },\text{ }\alpha >0\text{ ,}
\end{equation*}

\begin{eqnarray}  \label{eqn:Pex}
\mathbb{P}_{\varepsilon }\left( x\right) \leq \exp \left( -\frac{nx^{2}}{%
8\sigma_{\varepsilon }^{2}\left( 1+\alpha \right) }\right) + \mathbb{P}%
\left( \frac{1}{n}\left| \sum_{i=1}^{n}\left| \psi_{jk;s}\right|
^{2}-\left\| \psi _{jk;s}\right\| _{2}^{2}\right| >\alpha\right) \text{ .}
\notag
\end{eqnarray}

Now, by fixing $x=\kappa t_n$, we obtain the following result:
\begin{equation*}
\mathbb{P}_{\varepsilon }\left( \kappa t_n \right) \leq \exp \left( -\frac{%
k^{2}\log n}{8\sigma_{\varepsilon }^{2}\left( \alpha +1\right) }\right) +%
\mathbb{P}\left( \left| \frac{1}{n}\sum_{i=1}^{n}\left| \psi _{jk;s}\left(
X_{i}^{\prime }\right) \right| {^{2}-}\left\| \psi _{jk;s}\right\| {%
_{L_{s}^{2}\left( \mathbb{S}^{2}\right) }^{2}}\right| \geq \alpha \right)
\end{equation*}%
Now, we use on the second term the Hoeffding's inequality:
\begin{equation*}
\mathbb{P}\left( \left| \frac{1}{n}\sum_{i=1}^{n}\left| \psi _{jk;s}\left(
X_{i}^{\prime }\right) \right| ^{2}{-}\left\| \psi _{jk;s}\right\|
_{L_{s}^{2}\left( \mathbb{S}^{2}\right) }^{2}\right| \geq \alpha \right)
\leq 2\exp \left\{ -\frac{2n^{2}\alpha ^{2}}{ncB^{2j}}\right\} \text{ .}
\end{equation*}%
Again, because $B^{2j}\leq \frac{n}{\log n}$, we obtain:
\begin{eqnarray}  \label{eqn:Petn}
\mathbb{P}_{\varepsilon }\left( \kappa t_n\right) &\leq & 2\left\{ \exp
\left( -\frac{2\alpha ^{2}\log n}{c}\right) +\exp \left( -\frac{k^{2}\log n}{%
8\sigma _{\varepsilon }^{2}\left( \alpha +1\right) }\right) \right\} \\
&=&2\left\{ n^{-\frac{2\alpha ^{2}}{c}}+n^{-\frac{k^{2}}{8\sigma
_{\varepsilon }^{2}\left( \alpha +1\right) }}\right\} \text{ .}  \notag
\end{eqnarray}%
We fix $\alpha \sim k^{\frac{2}{3}}$ in order to obtain the same order of
magnitude between the two terms, and by using \ref{eqn:PFtn} and \ref%
{eqn:Petn} finally we obtain:
\begin{equation*}
\mathbb{P}_{\varepsilon },\mathbb{P}_{F}\leq C\cdot n^{-ck^{4/3}}\text{ .}
\end{equation*}

In order to prove \ref{eqn:prop2}, we use again \ref{eqn:split}, to obtain:%
\begin{equation*}
E\left[ \left| \widehat{\beta }_{jk;s}-\beta _{jk;s}\right| ^{p}\right]
\end{equation*}%
\begin{equation*}
\leq 2^{p-1}\left( E\left[ \left| \frac{1}{n}\sum_{i=1}^{n}\left( \overline{%
\psi }_{jk;s}F_{s}\left( X_{i}\right) \right) -E\left( \overline{\psi }%
_{jk;s}F_{s}\left( X_{i}\right) \right) \right| ^{p}\right] +E\left[ \left|
\frac{1}{n}\sum_{i=1}^{n}\overline{\psi }_{jk;s}\varepsilon _{i;s}\right|
^{p}\right] \right)
\end{equation*}%
\begin{equation*}
=2^{p-1}\left( E_{F}+E_{\varepsilon }\right) \text{ .}
\end{equation*}%
We need to split again both $E_{F}$ and $E_{\varepsilon }$ into real and
imaginary parts. Note that%
\begin{eqnarray*}
E_{F} &=&E\left[ \left| \frac{1}{n}\sum_{i=1}^{n}\func{Re}\Psi
_{jk;s}(X_{i})+\func{Im}\Psi _{jk;s}(X_{i})\right| ^{p}\right] \\
&\leq &2^{p-1}\left( E\left( \left| \frac{1}{n}\sum_{i=1}^{n}\func{Re}\Psi
_{jk;s}(X_{i})\right| ^{p}\right) +E\left( \left| \frac{1}{n}\sum_{i=1}^{n}%
\func{Im}\Psi _{jk;s}(X_{i})\right| ^{p}\right) \right) \\
&\leq &2^{p-1}\left( E_{F}^{1}+E_{F}^{2}\right)
\end{eqnarray*}%
and
\begin{eqnarray}
E_{\varepsilon } &=&E\left[ \left| \frac{1}{n}\sum_{i=1}^{n}\left( \func{Re}%
\left\{ \overline{\psi }_{jk;s}\left( X_{i}\right) \varepsilon
_{i;s}\right\} +\func{Im}\left\{ \overline{\psi }_{jk;s}\left( X_{i}\right)
\varepsilon _{i;s}\right\} \right) \right| ^{p}\right]  \notag \\
&\leq &2^{p-1}\left( E\left( \left| \frac{1}{n}\sum_{i=1}^{n}\func{Re}%
\left\{ \overline{\psi }_{jk;s}\left( X_{i}\right) \varepsilon
_{i;s}\right\} \right| ^{p}\right) +E\left( \left| \frac{1}{n}\sum_{i=1}^{n}%
\func{Im}\left\{ \overline{\psi }_{jk;s}\left( X_{i}\right) \varepsilon
_{i;s}\right\} \right| ^{p}\right) \right)  \label{eqn:Esplit} \\
&\leq &2^{p-1}\left( E_{\varepsilon }^{1}+E_{\varepsilon }^{2}\right) \text{
.}  \notag
\end{eqnarray}%
For $0<p\leq 2$, we apply the classical convexity inequality, which states
that for $0<p\leq 2$ a for independent random variables $Z_{i}$ such that $%
E\left( Z_{i}\right) =0$ and $E\left( \left| Z_{i}\right| ^{p}\right)
<\infty $:
\begin{equation*}
E\left( \left| \sum_{i=1}^{n}Z_{i}\right| ^{p}\right) \leq \left( E\left(
\left| \sum_{i=1}^{n}Z_{i}\right| ^{2}\right) \right) ^{\frac{p}{2}}\text{ .}
\end{equation*}%
As noted for instance in \cite{WASA}, in the case $2<p<\infty $, we obtain a
very similar result by applying the Rosenthal inequality, i.e.:

Let $\left\{ Z_{i}\right\} _{i=1}^{n}$ be independent random variable such
that $E\left( Z_{i}\right) =0$ and for $p\geq 2,$ $E\left( \left|
Z_{i}\right| ^{p}\right) <\infty $. Then there exists $C_{p}$ such that:
\begin{equation}
E\left( \left| \sum_{i=1}^{n}Z_{i}\right| ^{p}\right) \leq C_{p}\left[
\sum_{i=1}^{n}E\left( \left| Z_{i}\right| ^{p}\right) +\left(
\sum_{i=1}^{n}E\left( Z_{i}^{2}\right) \right) ^{\frac{p}{2}}\right] \text{ .%
}  \label{eqn:ros}
\end{equation}%
A proof of proposition (\ref{eqn:ros}) can be found for instance in \cite%
{WASA}. We apply (\ref{eqn:ros}) to each term in (\ref{eqn:Esplit}) to
obtain:
\begin{eqnarray}
E_{F}^{1} &\leq &C_{p}\left( \frac{E\left( \left| \func{Re}\Psi
_{jk;s}(X_{i})\right| ^{p}\right) }{n^{p-1}}+\frac{\left( E\left( \left|
\func{Re}\Psi _{jk;s}(X_{i})\right| ^{2}\right) \right) ^{\frac{p}{2}}}{n^{%
\frac{p}{2}}}\right) \text{ ,}  \label{eqn:rosf1} \\
E_{F}^{2} &\leq &C_{p}\left( \frac{E\left( \left| \func{Im}\Psi
_{jk;s}(X_{i})\right| ^{p}\right) }{n^{p-1}}+\frac{\left( E\left( \left|
\func{Im}\Psi _{jk;s}(X_{i})\right| ^{2}\right) \right) ^{\frac{p}{2}}}{n^{%
\frac{p}{2}}}\right) \text{ ,}  \label{eqn:rosf2} \\
E_{\varepsilon }^{1} &\leq &C_{p}\left( \frac{E\left( \left| \func{Re}\left(
\overline{\psi }_{jk;s}\left( X_{i}\right) \varepsilon _{i;s}\right) \right|
^{p}\right) }{n^{p-1}}+\frac{\left( E\left( \left| \func{Re}\left( \overline{%
\psi }_{jk;s}\left( X_{i}\right) \varepsilon _{i;s}\right) \right|
^{2}\right) \right) ^{\frac{p}{2}}}{n^{\frac{p}{2}}}\right) \text{ ,}
\label{eqn:rose1} \\
E_{\varepsilon }^{2} &\leq &C_{p}\left( \frac{E\left( \left| \func{Im}\left(
\overline{\psi }_{jk;s}\left( X_{i}\right) \varepsilon _{i;s}\right) \right|
^{p}\right) }{n^{p-1}}+\frac{\left( E\left( \left| \func{Im}\left( \overline{%
\psi }_{jk;s}\left( X_{i}\right) \varepsilon _{i;s}\right) \right|
^{2}\right) \right) ^{\frac{p}{2}}}{n^{\frac{p}{2}}}\right) \text{ .}
\label{eqn:rose2}
\end{eqnarray}%
Recalling that $B^{j}\leq \sqrt{\frac{n}{\log n}}\leq \sqrt{n}$, we obtain:
\begin{eqnarray*}
&&E\left( \func{Re}\left| \Psi _{jk;s}(X_{i})\right| ^{p}\right) =E\left(
\func{Im}\left| \Psi _{jk;s}(X_{i})\right| ^{p}\right) \leq \\
&\leq &E\left( \left| \overline{\psi }_{jk;s}\left( X_{i}\right) F_{s}\left(
X_{i}\right) \right| ^{p}\right) \leq \int_{\mathbb{S}^{2}}\left| \overline{%
\psi }_{jk;s}\left( X_{i}\right) F_{s}\left( X_{i}\right) \right| ^{p}dx \\
&\leq &cM^{p}B^{j\left( p-2\right) }\leq cM^{p}n^{-\frac{p-2}{2}}\text{ .}
\end{eqnarray*}%
As far as the noise-related terms, we obtain:
\begin{eqnarray*}
&&E\left( \left| \func{Re}\left( \overline{\psi }_{jk;s}\left( X_{i}\right)
\varepsilon _{i;s}\right) \right| ^{p}\right) =E\left( \left| \func{Im}%
\left( \overline{\psi }_{jk;s}\left( X_{i}\right) \varepsilon _{i;s}\right)
\right| ^{p}\right) \\
&\leq &E\left( \left| \varepsilon _{i;s}\right| ^{p}\right) cB^{j\left(
p-2\right) }\leq cn^{-\frac{p-2}{2}}\text{ .}
\end{eqnarray*}%
Then, by substituting the last inequalities in \ref{eqn:rosf1}, \ref%
{eqn:rosf2}, \ref{eqn:rose1} and \ref{eqn:rose2}, we obtain:
\begin{equation*}
\frac{n^{\frac{p}{2}-1}}{n^{p-1}}=n^{-\frac{p}{2}}
\end{equation*}%
Now we study the case $p=\infty $: in order to prove (\ref{eqn:prop3}), we
majorize:
\begin{equation}
E\left[ \sup_{k}\left| \widehat{\beta }_{jk;s}-\beta _{jk;s}\right| \right]
\leq \int_{\mathbb{R}^{+}}x^{p-1}\mathbb{P}\left( \sup_{k}\left| \widehat{%
\beta }_{jk;s}-\beta _{jk;s}\right| ^{p}>x\right) dx\text{ .}
\label{eqn:majsup}
\end{equation}

Recalling the procedure used in the proof of \ref{eqn:prop1}, for $B^{j}\leq
\sqrt{n}$, (\ref{eqn:PFx}) becomes:
\begin{equation}
\mathbb{P}_{F}\left( x\right) \leq 4\left( \exp \left( -\frac{nx^{2}}{16M^{2}%
}\right) +\exp \left( -\frac{3\sqrt{n}x}{16cM}\right) \right) \text{ ,}
\label{pf}
\end{equation}%
while, in a similar way, we split the first term on \ref{eqn:Pex} as:
\begin{equation}
\exp \left( -\frac{nx^{2}}{8\sigma _{\varepsilon }^{2}\left( 1+\alpha
\right) }\right) \leq \exp \left( -\frac{nx^{2}}{16\sigma _{\varepsilon }^{2}%
}\right) +\exp \left( -\frac{nx^{2}}{16\sigma _{\varepsilon }^{2}\alpha }%
\right) =\mathbb{P}_{\varepsilon }^{\ast }\left( x\right) +\mathbb{P}%
_{\varepsilon ,\alpha }^{1}\text{ .}  \label{pe3}
\end{equation}%
By applying on the last term of \ref{eqn:Pex} the Hoeffding inequality and
for $B^{j}\leq \sqrt{n}$, we obtain:
\begin{eqnarray*}
\mathbb{P}\left( \frac{1}{n}\left| \sum_{i=1}^{n}\left| \psi _{jk;s}\right|
^{2}-\left\| \psi _{jk;s}\right\| _{2}^{2}\right| >\alpha \right) &\leq
&\exp \left( -\frac{2n^{2}\alpha ^{2}}{ncB^{2j}}\right) \\
&\leq &\exp \left( -\frac{2\alpha ^{2}}{c}\right) =\mathbb{P}_{\varepsilon
,\alpha }^{2}\text{ .}
\end{eqnarray*}%
We choose $\alpha =\left\{ \frac{c^{1/3}}{32^{1/3} \sigma _{\varepsilon
}^{2/3}}\cdot n^{1/3}x^{2/3}\right\} $, to obtain%
\begin{equation}
\mathbb{P}_{\varepsilon ,\alpha }^{1}+\mathbb{P}_{\varepsilon ,\alpha
}^{2}\leq C\exp \left( -\frac{n^{2/3}x^{4/3}}{2^{7/3}\sigma _{\varepsilon
}^{4/3}c^{1/3}}\right) \text{ ,}  \label{pe2}
\end{equation}%
and in view of (\ref{pe1}), (\ref{pf}), (\ref{pe2}), (\ref{pe3})
\begin{equation*}
\mathbb{P}_{\beta }\left( x\right) \leq C\left( \exp \left( -\frac{nx^{2}}{%
16\sigma _{\varepsilon }^{2}}\right) +\exp \left( -\frac{nx^{2}}{16M^{2}}%
\right) +\exp \left( -\frac{2\sqrt{n}x}{16cM}\right) +\exp \left( -\frac{%
n^{2/3}x^{4/3}}{2^{7/3}\sigma _{\varepsilon }^{4/3}c^{1/3}}\right) \right)
\text{ .}
\end{equation*}%
Now we fix a parameter $a=max\left( 4\sqrt{2}\sigma _{\varepsilon },4\sqrt{2}%
M,\frac{32}{3}cM,2^{11/4}\sigma _{\varepsilon }c^{1/4}\right) $. Write (\ref%
{eqn:majsup}) as:
\begin{eqnarray}
&&E\left[ \sup_{k}\left|\widehat{\beta }_{jk;s}-\beta _{jk;s}\right| ^{p} %
\right] \leq \int_{0\leq x\leq \frac{aj}{\sqrt{n}}}x^{p-1}dx +  \notag \\
&& +2c \int_{x>\frac{aj}{\sqrt{n}}}C x^{p-1} B^{j}\left[ \exp \left( \frac{
nx^{2}}{16\sigma _{\varepsilon }^{2}}\right) + \exp \left( -\frac{nx^{2}}{%
16M^{2}}\right) + \exp \left( -\frac{2\sqrt{n}x}{16cM}\right) + \exp \left(
-\left(\frac{nx^2}{2^7 \sigma _{\varepsilon }^2 c^{\frac{1}{2}}} \right)^{%
\frac{2}{3}} \right)\right]  \notag \\
&=&E_{\infty }^{1}+E_{\infty }^{2}+E_{\infty }^{3}\text{ .}
\label{eqn:intsplit}
\end{eqnarray}%
We observe that for each term depending on $\exp \left( -nx^{2}/C\right) $,
where $C=4\sqrt{2}\sigma _{\varepsilon },4\sqrt{2}M$, and for $x>aj/\sqrt{n}$%
, we have:
\begin{equation*}
B^{j}\exp \left( -\frac{nx^{2}}{C}\right) \leq \exp \left( -\frac{nx^{2}}{2C}%
-\frac{nx^{2}}{2C}+j\right) \leq \exp \left( -\frac{nx^{2}}{2C}\right) \text{
.}
\end{equation*}%
Similarly, we have for $x>aj/\sqrt{n}$:
\begin{equation*}
B^{j}\exp \left( -\frac{2\sqrt{n}x}{16cM}\right) \leq \exp \left( -\frac{2%
\sqrt{n}x}{32cM}\right) \text{ ,}
\end{equation*}%
and finally, again for $x>aj/\sqrt{n}$
\begin{equation*}
B^{j}\exp \left( -\frac{n^{2/3}x^{4/3}}{2^7/3\sigma _{\varepsilon }^{4/3}}%
\right) \leq \exp \left( -\frac{n^{2/3}x^{4/3}}{2^{10^3}\sigma _{\varepsilon
}^{4/3}}\right) \text{ .}
\end{equation*}%
Likewise, the integral $E_{\infty }^{1}$ is simply majorized by:
\begin{equation}
E_{\infty }^{1}\leq C\frac{1}{p}\left( \frac{j}{\sqrt{n}}\right) ^{p}\leq
C_{p}j^{p}n^{-p/2}\text{ .}  \label{eqn:E1}
\end{equation}%
As far as $E_{\infty }^{2}$ is concerned, by using a change of variable $u=%
\sqrt{n}x$ we obtain:
\begin{equation}
E_{\infty }^{2}\leq 2C\frac{1}{n^{-p/2}}\int_{u\geq aj}u^{p-1}\exp \left( -%
\frac{u^{4/3}}{2^{10/3}\sigma _{\varepsilon }^{4/3}c^{1/3}}\right) du\leq
C_p n^{-p/2}.  \label{eqn:E2}
\end{equation}%
A similar procedure is applied to $E_{\infty }^{3}$ by using the same change
of variable $u=\sqrt{n}x$ to obtain:
\begin{equation}
E_{\infty }^{3}\leq C^{\prime}_p n^{-p/2}\text{ .}  \label{eqn:E3}
\end{equation}%
Finally by substituting (\ref{eqn:E1}), (\ref{eqn:E2}) and (\ref{eqn:E3}) in
(\ref{eqn:intsplit}) we obtain the thesis.\hfill $\square $

\subsection{Proof of Theorem \ref{th:loss}}

As customary in this literature, the proof can be divided into different
cases, as follows.

\begin{itemize}
\item Regular zone, $p<\infty $
\end{itemize}

We start as usual from
\begin{eqnarray*}
E\left\| F_{s}^{\ast }-F_{s}\right\| _{L_{s}^{p}}^{p} &=&E\left\|
\sum_{j\leq J_{n}}\sum_{k}w_{jk}\widehat{\beta }_{jk;s}\psi
_{jk;s}-\sum_{j}\sum_{k}\beta _{jk;s}\psi _{jk;s}\right\| _{L_{s}^{p}\left(
\mathbb{S}^{2}\right) }^{p} \\
&=&E\left\| \sum_{j\leq J_{n}}\sum_{k}(w_{jk}\widehat{\beta }_{jk;s}-\beta
_{jk;s})\psi _{jk;s}+\sum_{j>J_{n}}\sum_{k}\beta _{jk;s}\psi _{jk;s}\right\|
_{L_{s}^{p}\left( \mathbb{S}^{2}\right) }^{p} \\
&\leq &E\left\| \sum_{j\leq J_{n}}\sum_{k}(w_{jk}\widehat{\beta }%
_{jk;s}-\beta _{jk;s})\psi _{jk;s}\right\| _{L_{s}^{p}\left( \mathbb{S}%
^{2}\right)}^{p}+\left\| \sum_{j>J_{n}}\sum_{k}\beta _{jk;s}\psi
_{jk;s}\right\| _{L_{s}^{p}\left( \mathbb{S}^{2}\right) }^{p} \\
&=&:I+II\text{ .}
\end{eqnarray*}%
For $p\leq \pi ,$ we have $\mathcal{B}_{\pi q;s}^{r}\subset \mathcal{B}%
_{pq;s}^{r}$, whence we can always take $\pi =p$ in this case; hence we
focus on $p\geq \pi .$ Here we have the embedding $\mathcal{B}_{\pi
q;s}^{r}\subset \mathcal{B}_{pq;s}^{r-\frac{2}{\pi }+\frac{2}{p}},$ whence%
\begin{equation*}
\left\| \sum_{j>J_{n}}\sum_{k}\beta _{jk;s}\psi _{jk;s}\right\|
_{L_{s}^{p}\left( \mathbb{S}^{2}\right) }=O\left( B^{-2j(\frac{r}{2}-\frac{1%
}{\pi }+\frac{1}{p})}\right) =O\left( \left\{ \frac{n}{\log n}\right\} ^{-(%
\frac{r}{2}-\frac{1}{\pi }+\frac{1}{p})}\right) \text{ ,}
\end{equation*}%
and because in the regular zone%
\begin{equation*}
r\geq \frac{2}{\pi }\text{ , }\frac{r}{2r+2}=\frac{rp}{2(r+1)p}\leq \frac{%
r\pi }{2p}\text{ ,}
\end{equation*}%
we obtain
\begin{equation*}
(\frac{r}{2}-\frac{1}{\pi }+\frac{1}{p})-\frac{r}{2r+2}\geq (\frac{r}{2}-%
\frac{1}{\pi }+\frac{1}{p})-\frac{r\pi }{2p}=(\frac{1}{\pi }-\frac{1}{p})(%
\frac{r\pi }{2}-1)>0\text{ .}
\end{equation*}%
Hence the bias term is fixed. For the variance term we have%
\begin{equation*}
I\leq J_{n}^{p-1}\sum_{j\leq J_{n}}E\left\| \sum_{k}(w_{jk}\widehat{\beta }%
_{jk;s}-\beta _{jk;s})\psi _{jk;s}\right\| _{L_{s}^{p}\left( \mathbb{S}%
^{2}\right) }^{p}
\end{equation*}%
Now we split $I$ in four zones; more precisely, we shall label $A$
(respectively $U$) where the estimated coefficients is above (resp. under)
the threshold $\kappa t_{n}$, and $a$ (respectively $u$) the regions where
the deterministic coefficients are above or under a new threshold, which is $%
\frac{\kappa }{2}t_{n}$ in $A$ and $2\kappa t_{n}$ in $U.$ We hence obtain%
\begin{equation*}
\sum_{j\leq J_{n}}E\left\| \sum_{k}(w_{jk}\widehat{\beta }_{jk;s}-\beta
_{jk;s})\psi _{jk;s}\right\| _{L_{s}^{p}\left( \mathbb{S}^{2}\right) }^{p}
\end{equation*}%
\begin{eqnarray*}
&=&\sum_{j\leq J_{n}}E\left\| \sum_{k}(w_{jk}\widehat{\beta }_{jk;s}-\beta
_{jk;s})\psi _{jk;s}\right\| _{L_{s}^{p}\left( \mathbb{S}^{2}\right) }^{p}%
\mathbb{I}_{\left\{ \left| \widehat{\beta }_{jk;s}\right| \geq \kappa
t_{n}\right\} }\mathbb{I}_{\left\{ \left| \beta _{jk;s}\right| \geq \kappa
t_{n}/2\right\} } \\
&&+\sum_{j\leq J_{n}}E\left\| \sum_{k}(w_{jk}\widehat{\beta }_{jk;s}-\beta
_{jk;s})\psi _{jk;s}\right\| _{L_{s}^{p}\left( \mathbb{S}^{2}\right) }^{p}%
\mathbb{I}_{\left\{ \left| \widehat{\beta }_{jk;s}\right| \geq \kappa
t_{n}\right\} }\mathbb{I}_{\left\{ \left| \beta _{jk;s}\right| \leq \kappa
t_{n}/2\right\} } \\
&&+\sum_{j\leq J_{n}}E\left\| \sum_{k}(w_{jk}\widehat{\beta }_{jk;s}-\beta
_{jk;s})\psi _{jk;s}\right\| _{L_{s}^{p}\left( \mathbb{S}^{2}\right) }^{p}%
\mathbb{I}_{\left\{ \left| \widehat{\beta }_{jk;s}\right| <\kappa
t_{n}\right\} }\mathbb{I}_{\left\{ \left| \beta _{jk;s}\right| \geq 2\kappa
t_{n}\right\} } \\
&&+\sum_{j\leq J_{n}}E\left\| \sum_{k}(w_{jk}\widehat{\beta }_{jk;s}-\beta
_{jk;s})\psi _{jk;s}\right\| _{L_{s}^{p}\left( \mathbb{S}^{2}\right) }^{p}%
\mathbb{I}_{\left\{ \left| \widehat{\beta }_{jk;s}\right| <\kappa
t_{n}\right\} }\mathbb{I}_{\left\{ \left| \beta _{jk;s}\right| \leq 2\kappa
t_{n}\right\} }\text{ .}
\end{eqnarray*}%
\begin{eqnarray*}
&\leq &C\left\{ \sum_{j\leq J_{n}}\sum_{k}\left\| \psi _{jk;s}\right\|
_{L_{s}^{p}\left( \mathbb{S}^{2}\right) }^{p}E\left[ \left| \widehat{\beta }%
_{jk;s}-\beta _{jk;s}\right| ^{p}\mathbb{I}_{\left\{ \left| \widehat{\beta }%
_{jk;s}\right| \geq \kappa t_{n}\right\} }\mathbb{I}_{\left\{ \left| \beta
_{jk;s}\right| \geq \kappa t_{n}/2\right\} }\right] \right. \\
&&+\sum_{j\leq J_{n}}\sum_{k}\left\| \psi _{jk;s}\right\| _{L_{s}^{p}\left(
\mathbb{S}^{2}\right) }^{p}E\left[ \left| \widehat{\beta }_{jk;s}-\beta
_{jk;s}\right| ^{p}\mathbb{I}_{\left\{ \left| \widehat{\beta }_{jk;s}\right|
\geq \kappa t_{n}\right\} }\mathbb{I}_{\left\{ \left| \beta _{jk;s}\right|
\leq \kappa t_{n}/2\right\} }\right] \\
&&+\sum_{j\leq J_{n}}\sum_{k}\left\| \psi _{jk;s}\right\| _{L_{s}^{p}\left(
\mathbb{S}^{2}\right) }^{p}\left| \beta _{jk;s}\right| ^{p}E\left[ \mathbb{I}%
_{\left\{ \left| \widehat{\beta }_{jk;s}\right| <\kappa t_{n}\right\} }%
\mathbb{I}_{\left\{ \left| \beta _{jk;s}\right| \geq 2\kappa t_{n}\right\} }%
\right] \\
&&\left. +\sum_{j\leq J_{n}}\sum_{k}\left\| \psi _{jk;s}\right\|
_{L_{s}^{p}\left( \mathbb{S}^{2}\right) }^{p}\left| \beta _{jk;s}\right|
^{p}E\left[ \mathbb{I}_{\left\{ \left| \widehat{\beta }_{jk;s}\right|
<\kappa t_{n}\right\} }\mathbb{I}_{\left\{ \left| \beta _{jk;s}\right| \leq
2\kappa t_{n}\right\} }\right] \right\}
\end{eqnarray*}%
\begin{equation*}
=Aa+Au+Ua+Uu\text{ .}
\end{equation*}%
This idea is the same as in \cite{bkmpAoSb}, where the regions are labelled
instead $Bb,Bs,Sb,Ss;$ we preferred to avoid $B$ and $b$ which have a
different use in the present work. Heuristically, the cross/terms $Au,Ua$
are easier to bound, as we can exploit quick decay of $\Pr \left\{ \left|
\widehat{\beta }_{jk;s}-\beta _{jk;s}\right| >\frac{1}{2}t_{n}\right\} $;
for $Aa,Uu$ the crucial bounds will be derived by the tail behaviour in the
Besov balls $\mathcal{B}_{pq;s}^{r}(G)$.

Note firstly that%
\begin{eqnarray*}
Aa &\leq &C\sum_{j\leq J_{n}}\sum_{k}\left\| \psi _{jk;s}\right\|
_{L_{s}^{p}\left( \mathbb{S}^{2}\right) }^{p}E\left| \widehat{\beta }%
_{jk;s}-\beta _{jk;s}\right| ^{p}\mathbb{I}_{\left\{ \left| \beta
_{jk;s}\right| \geq \kappa t_{n}/2\right\} } \\
&\leq &C\sum_{j\leq J_{n}}\sum_{k}B^{j(p-2)}\mathbb{I}_{\left\{\left| \beta
_{jk;s}\right| \geq \kappa t_{n}/2\right\}}E\left| \widehat{\beta }%
_{jk;s}-\beta _{jk;s}\right| ^{p}\text{ ;}
\end{eqnarray*}%
now from \ref{eqn:prop2} and \ref{eqn:pnorm} we know that%
\begin{equation*}
E\left| \widehat{\beta }_{jk;s}-\beta _{jk;s}\right| ^{p}\leq C_{p}n^{-p/2}%
\text{ , }\sum_{k}B^{j(p-2)}=O(B^{jp})\text{ .}
\end{equation*}%
Write
\begin{eqnarray*}
&&\sum_{j\leq J_{n}}\sum_{k}B^{j(p-2)}\mathbb{I}_{\left\{\left|
\beta_{jk;s}\right| \geq \kappa t_{n}/2\right\}}E\left| \widehat{\beta }%
_{jk;s}-\beta_{jk;s}\right| ^{p} \\
& \leq & C\left\{ n^{-p/2}\sum_{j\leq J_{1n}}\sum_{k}B^{j(p-2)}\mathbb{I}%
_{\left\{\left| \beta_{jk;s}\right| \geq \kappa t_{n}/2\right\}} + n^{-p/2}
\sum_{j>J_{1n}}\sum_{k}B^{j(p-2)} \mathbb{I}_{\left\{\left|
\beta_{jk;s}\right| \geq \kappa t_{n}/2\right\}} \right\} \\
&\leq & C \left\{
n^{-p/2}B^{pJ_{1n}}+n^{-p/2}\sum_{j>J_{1n}}\sum_{k}B^{j(p-2)} \mathbb{I}%
_{\left\{\left| \beta_{jk;s}\right| \geq \kappa t_{n}/2\right\}} \right\}
\text{ .}
\end{eqnarray*}%
Fix
\begin{equation*}
B^{J_{1n}}=\kappa ^{\prime }\left\{ \frac{n}{\log n}\right\} ^{\frac{1}{%
2(r+1)}},
\end{equation*}%
and note that we have
\begin{eqnarray*}
\sum_{j>J_{1n}}\sum_{k}B^{j(p-2)}\mathbb{I}_{\left\{\left|
\beta_{jk;s}\right| \geq \kappa t_{n}/2\right\}} &\leq
&\sum_{j>J_{1n}}\sum_{k}\left| \beta _{jk;s}\right| ^{p}B^{j(p-2)}\left\{
\kappa t_{n}/2\right\} ^{-p} \\
&\leq &\left\{ \frac{n}{\log n}\right\} ^{p/2}\sum_{j>J_{1n}}\left\{
\sum_{k}\left| \beta _{jk;s}\right| ^{p}\left\| \psi _{jk;s}\right\|
_{L_{s}^{p}\left( \mathbb{S}^{2}\right) }^{p}\right\} \text{ ,}
\end{eqnarray*}%
where
\begin{equation*}
\left\{ \sum_{k}\left| \beta _{jk;s}\right| ^{p}\left\| \psi _{jk;s}\right\|
_{L_{s}^{p}\left( \mathbb{S}^{2}\right) }^{p}\right\} \leq CB^{-prj},
\end{equation*}%
because by assumption $F_{s}\in \mathcal{B}_{pq;s}^{r}$ .Hence%
\begin{eqnarray*}
&&\left\{ \frac{n}{\log n}\right\} ^{p/2}\sum_{j>J_{1n}}\left\{
\sum_{k}\left| \beta _{jk;s}\right| ^{p}\left\| \psi _{jk;s}\right\|
_{L_{s}^{p}\left( \mathbb{S}^{2}\right) }^{p}\right\} \\
&\leq &\left\{ \frac{n}{\log n}\right\} ^{p/2}B^{-prJ_{1n}}\leq C\left\{
\frac{n}{\log n}\right\} ^{p/2}\left\{ \frac{n}{\log n}\right\} ^{-\frac{pr}{%
2(r+1)}} \\
&\leq &C\left\{ \frac{n}{\log n}\right\} ^{\frac{p(r+1)-pr}{2(r+1)}}\leq
C\left\{ \frac{n}{\log n}\right\} ^{\frac{p}{2(r+1)}}\leq CB^{pJ_{1n}}\text{
,}
\end{eqnarray*}%
and%
\begin{eqnarray*}
&&C\sum_{j\leq J_{n}}\sum_{k}B^{j(p-2)}\mathbb{I}_{\left\{\left|
\beta_{jk;s}\right| \geq \kappa t_{n}/2\right\}} E\left| \widehat{\beta }%
_{jk;s}-\beta _{jk;s}\right| ^{p} \\
&\leq &Cn^{-p/2}B^{pJ_{1n}}\leq C\left\{ \frac{n}{\log n}\right\} ^{\frac{p}{%
2(r+1)}}n^{-p/2}\leq C\left\{ \frac{n}{\log n}\right\} ^{\frac{-pr}{2(r+1)}}
\text{ .}
\end{eqnarray*}%
Hence the term $Aa$ is fixed. For the term $Uu,$ it suffices to observe that
\begin{eqnarray*}
Uu &\leq &C\sum_{j\leq J_{n}}\sum_{k}\left\| \psi _{jk;s}\right\|
_{L_{s}^{p}\left( \mathbb{S}^{2}\right) }^{p}\left| \beta _{jk;s}\right| ^{p}%
\mathbb{I}_{\left\{ \left| \beta _{jk;s}\right| \leq 2\kappa t_{n}\right\} }
\\
&\leq &C\left\{ \sum_{j\leq J_{1n}}\sum_{k}B^{j(p-2)}\left| 2\kappa
t_{n}\right| ^{p}\mathbb{+}\sum_{j>J_{1n}}\sum_{k}B^{j(p-2)}\left| \beta
_{jk;s}\right| ^{p}\right\} \\
&\leq &C\left\{ B^{pJ_{1n}}\left\{ \frac{n}{\log n}\right\} ^{-p/2}\mathbb{+}%
B^{-prJ_{1n}}\right\} \\
&\leq &C\left\{ \left[ \frac{n}{\log n}\right] ^{\frac{p}{2(r+1)}}\left[
\frac{n}{\log n}\right] ^{-\frac{p}{2}}\mathbb{+}\left[ \frac{n}{\log n}%
\right] ^{-\frac{pr}{2(r+1)}}\right\} =O\left( \left[ \frac{n}{\log n}\right]
^{-\frac{pr}{2(r+1)}}\right) \text{ .}
\end{eqnarray*}%
Now note that%
\begin{eqnarray*}
Au &\leq &C\sum_{j\leq J_{n}}\sum_{k}B^{j(p-2)}E\left[ \left| \widehat{\beta
}_{jk;s}-\beta _{jk;s}\right| ^{p}\mathbb{I}_{\left\{ \left| \widehat{\beta }%
_{jk;s}-\beta _{jk;s}\right| \geq \kappa t_{n}/2\right\} }\right] \\
&\leq &\sum_{j\leq J_{n}}\sum_{k}B^{j(p-2)}\left\{ E\left[ \left| \widehat{%
\beta }_{jk;s}-\beta _{jk;s}\right| ^{2p}\right] \right\} ^{1/2}\left\{
\mathbb{P}\left[ \left| \widehat{\beta }_{jk;s}-\beta _{jk;s}\right| \geq
\kappa t_{n}/2\right] \right\} ^{1/2}
\end{eqnarray*}%
and using (\ref{eqn:prop2})%
\begin{equation*}
Au\leq C n^{-p/2}B^{pJ_{n}}n^{-\gamma /2}\leq C n^{-p/2}\left[ \frac{n}{\log
n}\right] ^{p/2}n^{-\gamma /2}=C\left( \log n\right) ^{-\frac{p}{2}%
}n^{-\gamma /2}.
\end{equation*}%
Finally%
\begin{equation*}
Ua\leq \sum_{j\leq J_{n}}\sum_{k}\left\| \psi _{jk;s}\right\|
_{L_{s}^{p}\left( \mathbb{S}^{2}\right) }^{p}\left| \beta _{jk;s}\right|
^{p}E\left[ \mathbb{I}_{\left\{ \left| \widehat{\beta }_{jk;s}-\beta
_{jk;s}\right| >\kappa t_{n}\right\} }\right] \leq Cn^{-\gamma }\left\|
F_{s}\right\| _{L_{s}^{p}\left(\mathbb{S}^2 \right)}^{p}.
\end{equation*}%
Because obviously $n^{-\gamma }\leq n^{-\gamma /2}$ we have to choose $%
\gamma $ such that:

\begin{equation*}
n^{-\gamma /2}\leq n^{-\frac{pr}{2r+2}}\longrightarrow \gamma \geq \frac{pr}{%
r+1}\text{ .}
\end{equation*}%
We can hence take $\kappa \sim \gamma ^{3/4}$, which yields
\begin{equation*}
\kappa \geq C\left( \frac{pr}{r+1}\right) ^{\frac{3}{4}}.
\end{equation*}

\begin{itemize}
\item The case $p=\infty $
\end{itemize}

Assume first that $F_{s}\in \mathcal{B}_{\infty ,\infty ;s}^{r}.$ Then%
\begin{eqnarray*}
E\left\| F_{s}^{\ast }-F_{s}\right\| _{L_{s}^{\infty }\left( \mathbb{S}%
^{2}\right) } &\leq &E\left\| \sum_{j\leq J_{n}}\sum_{k}(w_{jk}\widehat{%
\beta }_{jk;s}-\beta _{jk;s})\psi _{jk;s}\right\| _{L_{s}^{\infty }\left(
\mathbb{S}^{2}\right) }+\left\| \sum_{j>J_{n}}\sum_{k}\beta _{jk;s}\psi
_{jk;s}\right\| _{L_{s}^{\infty }\left( \mathbb{S}^{2}\right) } \\
&=&I+II\text{ .}
\end{eqnarray*}%
For $II,$\ it is sufficient to note that
\begin{eqnarray*}
\left\| \sum_{j>J_{n}}\sum_{k}\beta _{jk;s}\psi _{jk;s}\right\|
_{L_{s}^{\infty }\left( \mathbb{S}^{2}\right) } &\leq &\sum_{j>J_{n}}\left\|
\sum_{k}\beta _{jk;s}\psi _{jk;s}\right\| _{L_{s}^{\infty }\left( \mathbb{S}%
^{2}\right) }=O\left( B^{-rJ_{n}}\right) \\
&=&O\left( \left[ \frac{n}{\log n}\right] ^{-r/2}\right) =O\left( \left[
\frac{n}{\log n}\right] ^{-r/2(r+1)}\right) \text{ .}
\end{eqnarray*}%
On the other hand,
\begin{eqnarray*}
E\left\| \sum_{j\leq J_{n}}\sum_{k}(w_{jk}\widehat{\beta }_{jk;s}-\beta
_{jk;s})\psi _{jk;s}\right\| _{L_{s}^{\infty }\left( \mathbb{S}^{2}\right) }
&\leq &\sum_{j\leq J_{n}}E\left\| \sum_{k}(w_{jk}\widehat{\beta }%
_{jk;s}-\beta _{jk;s})\psi _{jk;s}\right\| _{L_{s}^{\infty }\left( \mathbb{S}%
^{2}\right) } \\
&\leq &C\sum_{j\leq J_{n}}B^{j}E\left[ \sup_{k}\left| w_{jk}\widehat{\beta }%
_{jk;s}-\beta _{jk;s}\right| \right]
\end{eqnarray*}%
\begin{eqnarray*}
&\leq &C\sum_{j\leq J_{n}}B^{j}E\left[ \sup_{k}\left| \widehat{\beta }%
_{jk;s}-\beta _{jk;s}\right| \right] \mathbb{I}_{\left\{ \left| \beta
_{jk;s}\right| \geq \kappa t_{n}/2\right\} } \\
&&+C\sum_{j\leq J_{n}}B^{j}E\left[ \sup_{k}\left| \widehat{\beta }%
_{jk;s}-\beta _{jk;s}\right| \right] \mathbb{I}_{\left\{ \left| \widehat{%
\beta }_{jk;s}-\beta _{jk;s}\right| \geq \kappa t_{n}/2\right\} } \\
&&+C\sum_{j\leq J_{n}}B^{j}\sup_{k}\left| \beta _{jk;s}\right| E\left[
\mathbb{I}_{\left\{ \left| \widehat{\beta }_{jk;s}-\beta _{jk;s}\right|
>\kappa t_{n}\right\} }\right] \\
&&+C\sum_{j\leq J_{n}}B^{j}\sup_{k}\left| \beta _{jk;s}\right| \mathbb{I}%
_{\left\{ \left| \beta _{jk;s}\right| \leq 2\kappa t_{n}\right\} } \\
&=&Aa+Au+Ua+Uu\text{ .}
\end{eqnarray*}%
Now as before, we note that it is possible to choose
\begin{equation*}
J_{1n}:B^{J_{1n}}\sim \kappa ^{\prime }\left\{ \frac{n}{\log n}\right\} ^{%
\frac{1}{2(r+1)}}\text{ and for }j>J_{1n}\text{ , }\mathbb{I}_{\left\{
\left| \beta _{jk;s}\right| \geq \kappa t_{n}/2\right\} }\equiv 0\text{ .}
\end{equation*}%
Hence, by (\ref{eqn:prop3})
\begin{eqnarray*}
Aa &\leq &C\sum_{j\leq J_{n}}B^{j}E\left[ \sup_{k}\left| \widehat{\beta }%
_{jk;s}-\beta _{jk;s}\right| \right] \mathbb{I}_{\left\{ \left| \beta
_{jk;s}\right| \geq \kappa t_{n}/2\right\} } \\
&\leq &C\sum_{j\leq J_{1n}}B^{j}E\left[ \sup_{k}\left| \widehat{\beta }%
_{jk;s}-\beta _{jk;s}\right| \right] \leq CJ_{1n}n^{-\frac{1}{2}}B^{J_{1n}}
\\
&\leq &CJ_{1n}\left( \log n\right) ^{-1/2}\left\{ \frac{n}{\log n}\right\}
^{-\frac{r}{2(r+1)}}.
\end{eqnarray*}

Also%
\begin{eqnarray*}
\sum_{j\leq J_{n}}B^{j}\sup_{k}\left| \beta _{jk}\right| \mathbb{I}_{\left\{
\left| \beta _{jk;s}\right| \leq 2\kappa t_{n}\right\} } &\leq &C\left\{
t_{n}B^{J_{1n}}+\sum_{J_{1n}\leq j<\infty }B^{j}\sup_{k}\left| \beta
_{jk}\right| \right\} \\
&\leq &C\left\{ t_{n}B^{J_{1n}}+\sum_{J_{1n}\leq j<\infty }\left\|
F_{s}\right\| _{L_{s}^{\infty }}\right\} \\
&\leq &C\left\{ t_{n}B^{J_{1n}}+B^{-J_{1n}}\right\} \leq C\left\{ \frac{n}{%
\log n}\right\} ^{-\frac{r}{2(r+1)}}.
\end{eqnarray*}%
For the remaining two terms the arguments is the same, actually easier. For
general $\pi $ and $q,$ it is sufficient to note that $\mathcal{B}_{\pi
q;s}^{r}\subset \mathcal{B}_{\infty ,\infty ;s}^{r^{\prime }},$ $r^{\prime
}=r-2/\pi .$ By the previous argument%
\begin{equation*}
E\left\| F_{s}^{\ast }-F_{s}\right\| _{L_{s}^{\infty }}\leq CJ_{n}\left\{
\frac{n}{\log n}\right\} ^{-\frac{r^{\prime }}{2(r^{\prime }+1)}%
}=CJ_{n}\left\{ \frac{n}{\log n}\right\} ^{-\frac{r-2/\pi }{2(r-2(1/\pi -1/2)%
}}.
\end{equation*}%
Note that for $\pi =p=\infty $ the sparse and regular zone coincide;
otherwise for $p=\infty $ we are always in the sparse zone

\begin{itemize}
\item The sparse case
\end{itemize}

The argument is very much the same as before. Indeed we have $\mathcal{B}%
_{\pi q;s}^{r}\subset \mathcal{B}_{p,q;s}^{r-2(\frac{1}{\pi }-\frac{1}{p})},$
\begin{equation*}
E\left\| F_{s}^{\ast }-F_{s}\right\| _{L_{s}^{p}}^{p}\leq E\left\|
\sum_{j\leq J_{n}}\sum_{k}(w_{jk}\widehat{\beta }_{jk;s}-\beta _{jk;s})\psi
_{jk;s}\right\| _{L_{s}^{p}\left( \mathbb{S}^{2}\right) }^{p}+\left\|
\sum_{j>J_{n}}\sum_{k}\beta _{jk;s}\psi _{jk;s}\right\| _{L_{s}^{p}\left(
\mathbb{S}^{2}\right) }^{p},
\end{equation*}%
\begin{eqnarray*}
\left\| \sum_{j>J_{n}}\sum_{k}\beta _{jk;s}\psi _{jk;s}\right\|
_{L_{s}^{p}\left( \mathbb{S}^{2}\right) }^{p} &\leq &CB^{-J_{n}(r-2(\frac{1}{%
\pi }-\frac{1}{p}))}\leq CB^{-2J_{n}[(r-2(\frac{1}{\pi }-\frac{1}{p}))/2(r-2(%
\frac{1}{\pi }-\frac{1}{p}))]} \\
&\leq &\left\{ \frac{n}{\log n}\right\} ^{-[(r-2(\frac{1}{\pi }-\frac{1}{p}%
))/2(r-2(\frac{1}{\pi }-\frac{1}{p}))]},
\end{eqnarray*}%
because $r-\frac{2}{\pi }+1\geq 1,$ given that $r-\frac{2}{\pi }\geq 0$ by
assumption. Hence the bias term has the correct order. For the variance
term, the trick is very much as above, and we omit some details. It is
possible to split the term to be bounded into four terms, after which the
two \ ''cross terms'' $Au$ and $Ua$ are easy because they involve quantities
like $\mathbb{P}\left\{ |\widehat{\beta }_{jk;s}-\beta _{jk;s}|>\kappa
t_{n}\right\} ,$ which can be made smaller than $n^{-p/2}$ for all $p>0,$
given a suitable choice of $\kappa .$ Fix $J_{2n}$ such that%
\begin{equation*}
B^{J_{2n}}\approx \left[ \frac{n}{\log n}\right] ^{\frac{1}{2((r-\frac{2}{%
\pi })+1)}},
\end{equation*}%
so that
\begin{eqnarray*}
\left[ \frac{n}{\log n}\right] ^{\frac{\pi -p}{2}}B^{J_{2n}(p-\pi (r+1))}
&\approx &\left[ \frac{n}{\log n}\right] ^{\frac{\pi -p}{2}}\left[ \frac{n}{%
\log n}\right] ^{\frac{(p-\pi (r+1))}{2((r-\frac{2}{\pi })+1)}} \\
&\approx &\left[ \frac{n}{\log n}\right] ^{\frac{(\pi -p)((r-\frac{2}{\pi }%
)+1)|+(p-\pi (r+1))}{2((r-\frac{2}{\pi })+1)}}.
\end{eqnarray*}%
For the terms of the form $Aa$ and $Uu\ $\ we have
\begin{equation*}
J_{n}^{p-1}n^{-p/2}\sum_{j\leq J_{1n}}B^{j(p-2)}\sum_{k}\mathbb{I}_{\left\{
\left| \beta _{jk;s}\right| \geq \kappa t_{n}/2\right\}
}+J_{n}^{p-1}\sum_{j}B^{j(p-2)}\sum_{k}|\beta _{jk;s}|^{p}\mathbb{I}%
_{\left\{ \left| \beta _{jk;s}\right| \leq 2\kappa t_{n}\right\} }\text{ ,}
\end{equation*}%
where to obtain the first summand we have exploited the embedding $\mathcal{B%
}_{\pi q;s}^{r}\subset \mathcal{B}_{\infty ,\infty ;s}^{r-\frac{2}{\pi }},$
whence for $j\geq J_{2n}$ one has $\mathbb{I}_{\left\{ \left| \beta
_{jk;s}\right| \geq \kappa t_{n}/2\right\} }\equiv 0.$ Now%
\begin{eqnarray*}
&&n^{-p/2}\sum_{j\leq J_{2n}}B^{j(p-2)}\sum_{k}\mathbb{I}_{\left\{ \left|
\beta _{jk;s}\right| \geq \kappa t_{n}/2\right\} } \\
&\leq &Cn^{-p/2}\sum_{j\leq J_{2n}}B^{j(p-2)}\sum_{k}|\beta _{jk;s}|^{\pi
}t_{n}^{-\pi } \\
&\leq &Cn^{-p/2}t_{n}^{-\pi }\sum_{j\leq J_{2n}}B^{j(p-\pi
)}\sum_{k}B^{j(\pi -2)}|\beta _{jk;s}|^{\pi } \\
&\leq &Cn^{-p/2}t_{n}^{-\pi }\sum_{j\leq J_{2n}}B^{j(p-\pi )}B^{-r\pi j}\leq
C\left[ \frac{n}{\log n}\right] ^{\frac{\pi -p}{2}}B^{J_{2n}(p-\pi (r+1))}.
\end{eqnarray*}%
Likewise%
\begin{eqnarray*}
\sum_{j\leq J_{2n}}B^{j(p-2)}\sum_{k}|\beta _{jk;s}|^{p}\mathbb{I}_{\left\{
\left| \beta _{jk;s}\right| \leq 2\kappa t_{n}\right\} } &\leq &C\sum_{j\leq
J_{2n}}B^{j(p-2)}\sum_{k}|\beta _{jk;s}|^{\pi }t_{n}^{p-\pi } \\
&\leq &C\left[ \frac{n}{\log n}\right] ^{\frac{\pi -p}{2}}\sum_{j\leq
J_{2n}}B^{j(p-\pi )}\sum_{k}B^{j(\pi -2)}|\beta _{jk;s}|^{\pi } \\
&\leq &C\left[ \frac{n}{\log n}\right] ^{\frac{\pi -p}{2}}B^{J_{2n}(p-\pi
(r+1))}.
\end{eqnarray*}%
Now%
\begin{eqnarray*}
\frac{(\pi -p)((r-\frac{2}{\pi })+1)+(p-\pi (r+1))}{2((r-\frac{2}{\pi })+1)}
&=&\frac{(\pi (r+1)-2-pr+\frac{2p}{\pi }-p)+(p-\pi (r+1))}{2((r-\frac{2}{\pi
})+1)} \\
&=&-\frac{2+pr-\frac{2p}{\pi }}{2((r-\frac{2}{\pi })+1)}=-\frac{p(r-2(\frac{1%
}{\pi }-\frac{1}{p}))}{2(r-2(\frac{1}{\pi }-\frac{1}{2}))}\text{ ,}
\end{eqnarray*}%
that is, these terms have the right order. So we are only left with
\begin{equation*}
\sum_{j\geq J_{2n}}B^{j(p-2)}\sum_{k}|\beta _{jk;s}|^{p}\mathbb{I}_{\left\{
\left| \beta _{jk;s}\right| \leq 2\kappa t_{n}\right\} }\text{ .}
\end{equation*}%
Consider%
\begin{equation*}
m=\frac{p-2}{r-\frac{2}{\pi }+1}\text{ ; }
\end{equation*}%
note that%
\begin{eqnarray*}
p-m &=&\frac{pr-\frac{2p}{\pi }+p-p+2}{r-\frac{2}{\pi }+1} \\
&=&\frac{pr-\frac{2p}{\pi }+2}{r-\frac{2}{\pi }+1}>0\text{ ,} \\
m-\pi  &=&\frac{p-2}{r-\frac{2}{\pi }+1}-\pi  \\
&=&\frac{p-\pi (r+1)}{r-\frac{2}{\pi }+1}>0\text{ ,}
\end{eqnarray*}%
because $p-\pi (r+1)>0$ in the sparse zone. We have%
\begin{eqnarray}
\sum_{j\geq J_{2n}}B^{j(p-2)}\sum_{k}|\beta _{jk;s}|^{p}\mathbb{I}_{\left\{
\left| \beta _{jk;s}\right| \leq 2\kappa t_{n}\right\} } &\leq &C\sum_{j\geq
J_{2n}}B^{j(p-2)}\sum_{k}|\beta _{jk;s}|^{m}t_{n}^{p-m}  \notag \\
&\leq &Ct_{n}^{p-m}\sum_{j\geq J_{2n}}B^{j(p-m)}\sum_{k}B^{j(m-2)}|\beta
_{jk;s}|^{m}  \notag \\
&\leq &Ct_{n}^{p-m}\sum_{j\geq J_{2n}}B^{j(p-m)}\sum_{k}\left\| \psi
_{jk;s}\right\| _{L_{s}^{m}\left( \mathbb{S}^{2}\right) }^{m}|\beta
_{jk;s}|^{m}.  \label{hot}
\end{eqnarray}%
Now, because $\mathcal{B}_{\pi q;s}^{r}\subset \mathcal{B}_{m,q;s}^{r-\frac{2%
}{\pi }+\frac{2}{m}},$%
\begin{equation*}
\sum_{k}\left\| \psi _{jk;s}\right\| _{L_{s}^{m}\left( \mathbb{S}^{2}\right)
}^{m}|\beta _{jk;s}|^{m}\leq CB^{-mj(r-\frac{2}{\pi }+\frac{2}{m})},
\end{equation*}%
hence (\ref{hot}) is bounded by%
\begin{equation*}
Ct_{n}^{p-m}\sum_{J_{2n}\leq j\leq J}B^{j(p-m-2)}B^{-j(r-\frac{2}{\pi }+%
\frac{2}{m})m}\leq Ct_{n}^{p-m}\sum_{J_{2n}\leq j\leq J}B^{j\left[
(p-m-2)-(r-\frac{2}{\pi }+\frac{2}{m})m\right] }.
\end{equation*}%
Observe that
\begin{eqnarray*}
(p-m)-(r-\frac{2}{\pi }+\frac{2}{m})m &=&\frac{pr-\frac{2p}{\pi }+2}{r-\frac{%
2}{\pi }+1}-(r-\frac{2}{\pi })m-2 \\
&=&\frac{pr-\frac{2p}{\pi }+2}{r-\frac{2}{\pi }+1}-(r-\frac{2}{\pi })\frac{%
p-2}{r-\frac{2}{\pi }+1}-2 \\
&=&\frac{2r+2(1-\frac{2}{\pi })}{r-\frac{2}{\pi }+1}-2=0\text{ ,}
\end{eqnarray*}%
hence%
\begin{equation*}
\sum_{J_{2n}\leq j\leq J}B^{j(p-2)}\sum_{k}|\beta _{jk;s}|^{p}\mathbb{I}%
_{\left\{ \left| \beta _{jk;s}\right| \leq 2\kappa t_{n}\right\} }\leq
CJ_{n}t_{n}^{p-m}\leq C\log n\left[ \frac{n}{\log n}\right] ^{-\frac{p(r-2(%
\frac{1}{\pi }-\frac{1}{p}))}{2(r-2(\frac{1}{\pi }-\frac{1}{2}))}}.
\end{equation*}%
Thus the proof is completed.\hfill $\square $

\end{document}